\documentclass[10pt]{amsart}

\usepackage{amsmath,amsthm,verbatim,amssymb,amsfonts,amscd, graphicx}
\usepackage{mathrsfs}
\usepackage[english]{babel}
\usepackage[retainorgcmds]{IEEEtrantools}
\usepackage{tikz-cd}
\usepackage{graphics}
\usepackage[all]{xy}
\usepackage{incgraph}
\usepackage{mathtools}
\usepackage{graphicx}
\usepackage{marginnote}
\usepackage{tikz}
\usetikzlibrary{arrows.meta}
\usetikzlibrary{bending}

\usepackage{geometry}
\usepackage{hyperref}

\newtheorem{thm}{Theorem}[section]

\newtheorem{thm-defn}[thm]{Theorem/Definition}
\newtheorem{lem}[thm]{Lemma}
\newtheorem{prop}[thm]{Proposition}
\newtheorem{cor}[thm]{Corollary}

\theoremstyle{definition}
\newtheorem{defn}[thm]{Definition}

\newtheorem{exam}[thm]{Example}

\newtheorem{set-up}[thm]{Set-up}

\theoremstyle{remark}
\newtheorem{rem}[thm]{Remark}

\numberwithin{equation}{section}

\def\Hom{\mathrm{Hom}}

\def\Cp{\mathbb{C}_p}

\def\E{\mathcal E} 

\def\makeop#1{\expandafter\def\csname#1\endcsname
  {\mathop{\rm #1}\nolimits}\ignorespaces}
\makeop{Hom}   \makeop{End}   \makeop{Aut}   \makeop{Isom}  \makeop{Pic}
\makeop{Gal}   \makeop{ord}   \makeop{Char}  \makeop{Div}   \makeop{Lie}
\makeop{PGL}   \makeop{Corr}  \makeop{PSL}   \makeop{sgn}   \makeop{Spf}
\makeop{Spec}  \makeop{Map}    \makeop{Nm}    \makeop{Fr}    \makeop{disc}
\makeop{Proj}  \makeop{supp}  \makeop{ker}   \makeop{im}    \makeop{dom}
\makeop{coker} \makeop{Stab}  \makeop{SO}    \makeop{SL}    \makeop{SL}
\makeop{Cl}    \makeop{cond}  \makeop{Br}    \makeop{inv}   \makeop{rank}
\makeop{id}    \makeop{Fil}   \makeop{Frac}  \makeop{GL}    \makeop{SU}
\makeop{Trd}   \makeop{Sp}    \makeop{Tr}    \makeop{Trd}   \makeop{diag}
\makeop{Res}   \makeop{ind}   \makeop{depth} \makeop{Tr}    \makeop{st}
\makeop{Ad}    \makeop{Int}   \makeop{tr}    \makeop{Sym}   \makeop{can}
\makeop{length}\makeop{SO}    \makeop{torsion} \makeop{GSp} \makeop{Ker}
\makeop{Adm}   \makeop{Frob}  \makeop{id}    \makeop{Tor}   \makeop{Ind}
\makeop{CoInd} \makeop{Inf}   \makeop{Cor}   \makeop{CoInf} \makeop{coLie}
\makeop{Ann}   \makeop{Shv}   \makeop{perf}  \makeop{Cont}  \makeop{Perfd}
\makeop{can}

\newcommand{\Z}{\mathbb Z}
\newcommand{\Q}{\mathbb Q}

\newcommand{\C}{\mathbb C}
\newcommand{\N}{\mathbb N}
\newcommand{\D}{\mathbf D}    
\renewcommand{\O}{\mathcal O} 

\newcommand{\BdR}{B_{\mathrm{dR}}} 
\newcommand{\Ainf}{A_{\mathrm{inf}}} 
\newcommand{\Acris}{A_{\mathrm{cris}}} 
\newcommand{\Bcris}{B_{\mathrm{cris}}} 
\newcommand{\Bst}{B_{\mathrm{st}}} 
\newcommand{\Spa}{\mathrm{Spa}} 
 
\newcommand{\MM}{\mathfrak{M}}
\newcommand{\SSS}{\mathfrak{S}}

\newcommand{\m}{\mathfrak m}

\newcommand{\p}{\mathfrak p}

\newcommand{\Spd}{\mathrm{Spd}}

\newcommand{\Fi}{\mathrm{Fil}}

\newcommand{\MMM}{\mathfrak{M}^\mathrm{inf}}

\newcommand{\cris}{\mathrm{cris}}
\newcommand{\lcr}{\mathrm{lcr}}
\newcommand{\dR}{\mathrm{dR}}

\newcommand{\MF}{\mathbf{MF}}
\newcommand{\BKF}{\mathbf{BKF}}
\newcommand{\HT}{\mathbf{HT}}
\newcommand{\Rep}{\mathrm{Rep}}
\newcommand{\cKur}{\breve{K}}

\usepackage{relsize}
\usepackage[bbgreekl]{mathbbol}
\usepackage{amsfonts}
\DeclareSymbolFontAlphabet{\mathbb}{AMSb} 
\DeclareSymbolFontAlphabet{\mathbbl}{bbold}
\newcommand{\Prism}{{\mathlarger{\mathbbl{\Delta}}}}

\newcommand{\M}{\mathcal M}

\newcommand{\Bun}{\mathrm{Bun}} 
\newcommand{\Modif}{\mathscr{M}odif} 
\newcommand{\phimod}{\varphi\text{-}\mathrm{Mod}}

\makeatletter
\newcommand{\colim@}[2]{%
  \vtop{\m@th\ialign{##\cr
    \hfil$#1\operator@font colim$\hfil\cr
    \noalign{\nointerlineskip\kern1.5\ex@}#2\cr
    \noalign{\nointerlineskip\kern-\ex@}\cr}}%
}
\newcommand{\colim}{%
  \mathop{\mathpalette\colim@{\rightarrowfill@\textstyle}}\nmlimits@
}
\makeatother

\makeatletter
\newcommand*{\da@rightarrow}{\mathchar"0\hexnumber@\symAMSa 4B }
\newcommand*{\da@leftarrow}{\mathchar"0\hexnumber@\symAMSa 4C }
\newcommand*{\xdashrightarrow}[2][]{%
  \mathrel{%
    \mathpalette{\da@xarrow{#1}{#2}{}\da@rightarrow{\,}{}}{}%
  }%
}
\newcommand{\xdashleftarrow}[2][]{%
  \mathrel{%
    \mathpalette{\da@xarrow{#1}{#2}\da@leftarrow{}{}{\,}}{}%
  }%
}
\newcommand*{\da@xarrow}[7]{%
  \sbox0{$\ifx#7\scriptstyle\scriptscriptstyle\else\scriptstyle\fi#5#1#6\m@th$}%
  \sbox2{$\ifx#7\scriptstyle\scriptscriptstyle\else\scriptstyle\fi#5#2#6\m@th$}%
  \sbox4{$#7\dabar@\m@th$}%
  \dimen@=\wd0 %
  \ifdim\wd2 >\dimen@
    \dimen@=\wd2 %
  \fi
  \count@=2 %
  \def\da@bars{\dabar@\dabar@}%
  \@whiledim\count@\wd4<\dimen@\do{%
    \advance\count@\@ne
    \expandafter\def\expandafter\da@bars\expandafter{%
      \da@bars
      \dabar@ 
    }%
  }%
  \mathrel{#3}%
  \mathrel{%
    \mathop{\da@bars}\limits
    \ifx\\#1\\%
    \else
      _{\copy0}%
    \fi
    \ifx\\#2\\%
    \else
      ^{\copy2}%
    \fi
  }%
  \mathrel{#4}%
}
\makeatother
\begin{document}
\author{Heng Du}
\address[Heng Du]{Yau Mathematical Sciences Center, Tsinghua University, Beijing 100084, China}
\email{hengdu@mail.tsinghua.edu.cn}

\title[Arithmetic Breuil-Kisin-Fargues modules and comparison of integral p-adic Hodge theories]{Arithmetic Breuil-Kisin-Fargues modules and comparison of integral p-adic Hodge theories}

\maketitle

\begin{abstract}
Let $K$ be a discrete valuation field with perfect residue field, we study the functor from weakly admissible filtered $(\varphi,N,G_K)$-modules over $K$ to the isogeny category of Breuil-Kisin-Fargues $G_K$-modules. This functor is the composition of a functor defined by Fargues-Fontaine from weakly admissible filtered $(\varphi,N,G_K)$-modules to $G_K$-equivariant modifications of vector bundles over the Fargues-Fontaine curve $X_{FF}$, with the functor of Fargues-Scholze that between the category of admissible modifications of vector bundles over $X_{FF}$ and the isogeny category of Breuil-Kisin-Fargues modules. We study those objects that appear in the essential image of the above functor and call them arithmetic BKF modules. We show certain rigidity result of arithmetic BKF modules and use it to compare existing $p$-adic Hodge theories at $\Ainf$ level.
\end{abstract}

\section{Introduction}

\subsection{Review of the work of Fargues-Fontaine and Fargues-Scholze} Fargues and Fontaine in \cite{FF} construct a complete abstract curve $X_{FF}$, the Fargues-Fontaine curve (constructed using the perfectoid field $\Cp^\flat$ and $p$-adic field $\Q_p$). For any $p$-adic field $K$, they show $\O_X=\O_{X_{FF}}$ carries an action of $G_K$, and they define $\O_X$-representations of $G_K$ as vector bundles over $X_{FF}$ that carries a continuous $\O_X$-semilinear action of $G_K$. They can show $\O_X$-representations of $G_K$ can help study $p$-adic representations of $G_K$ in many aspects. For example, Fargues-Fontaine show there is a nice slope theory on $X_{FF}$, and prove that the category of $\O_X$-representations such that the underlying vector bundles over $X_{FF}$ are semistable of slope 0 is equivalence to the category of $p$-adic Galois representations. Moreover, they give an explicit construction of slope 0 $\O_X$-representations from weakly admissible filtered $(\varphi,N)$-modules $D$ over $K$. Their construction is that: first using $D$ and the $(\varphi,N)$-structure, they construct an $\O_X$-representation $\E(D,\varphi,N)$ of $G_K$ whose underlying vector bundle is not semistable in general, then using the filtration structure of $D_K$, they constructed a $G_K$-equivariant modification $\E(D,\varphi,N,\Fi^\bullet)$ of $\E(D,\varphi,N)$ along the $G_K$-stable closed point $\infty$ on $X_{FF}$. They can show if $D$ is weakly admissible, then $\E(D,\varphi,N,\Fi^\bullet)$ is of slope $0$, and the $\Q_p$-representation corresponds to $\E(D,\varphi,N,\Fi^\bullet)$ is nothing but the log-crystalline representation corresponding to the data $(D,\varphi,N,\Fi^\bullet)$. Using this construction, they give new proofs of some important theorems in $p$-adic Hodge theory, for instance, the result that weakly admissible implies admissible, and also the $p$-adic monodromy theory of $p$-adic Galois representations.

One can also view the work of Fargues-Fontaine by introducing the following definitions.

\begin{defn}[Fargues]
The category of $p$-adic Hodge structures over $\C_p$ is defined as the category of admissible modifications of vector bundles over $X_{FF}$. 
\end{defn}

And we make the following definition.

\begin{defn}
\begin{enumerate}
    \item The category of $p$-adic Hodge structures over $\C_p$ with $G_K$-structures is defined as the category of admissible $G_K$-equivariant modifications of vector bundles over $X_{FF}$;
    \item The category of (arithmetic) $p$-adic Hodge structures over $K$ is just the category of weakly admissible filtered $(\varphi,N,G_K)$-modules over $K$.
\end{enumerate} 
\end{defn}

And Fargues-Fontaine in \cite{FF} constructed a functor $\eta_{FF}$ from the category of arithmetic $p$-adic Hodge structures over $K$ to the category of $p$-adic Hodge structures over $\C_p$ with $G_K$-structures, which we will call the Fargues-Fontaine functor. And the first result of this paper is that 

\begin{thm}[Theorem~\ref{etaFFisff}]\label{thm-intro-FFfunctor}
The Fargues-Fontaine functor $\eta_{FF}$ is fully faithful and the essential image of $\eta_{FF}$ consists of those admissible $G_K$-equivariant modifications $(\E_0 \xdashrightarrow[]{} \E_1)$ such that the fiber of $\E_1$ at $\infty$ as a $\C_p$-representation of $G_K$ is $\C_p$ admissible. 
\end{thm}

\subsection{Arithmetic Breuil-Kisin-Fargues modules and essential images of Fargues-Fontaine-Scholze functor} 
The abstract curve $X_{FF}$ also plays an important role in Scholze's work. In his Berkeley lecture on $p$-adic geometry\cite{ScholzeWeinstein}, Scholze defined a mixed characteristic analog of shtukas with legs. To be more precise, he introduced the functor $\Spd(\Z_p)$ which plays a similar role of a proper smooth curve in the equal characteristic story, and for any perfectoid space $S$ in characteristic $p$, he was able to define the category of shtukas over $S$. If we restrict us to the case that when $S=\Spa(C)$ with $C=\Cp^\flat$, and assume there is just one leg at the point corresponds to the untilt $\Cp$, then he can realize shtukas over $S$ as (admissible) modifications of vector bundles over $X_{FF}$ along $\infty$. Here $\infty$ is the same closed point on $X_{FF}$ as we mentioned in the work of Fargues-Fontaine. Fargues-Scholze also show that those shtukas can be realized using some commutative algebra data, called finite free Breuil-Kisin-Fargues modules, which are modules over $\Ainf=W(\O_{C})$ with some additional structures.

If we combine the construction of Fargues and Fontaine of modifications of vector bundles over $X_{FF}$ from log-crystalline representations and the work of Fargues and Scholze that relates modifications of vector bundles over $X_{FF}$ with local shtukas and Breuil-Kisin-Fargues modules, one can expect that if starting with (integral) arithmetic $p$-adic Hodge structures over $K$, one can produce a finite free Breuil-Kisin-Fargues module using the admissible modification constructed by Fargues-Fontaine. Moreover, since the modification is $G_K$-equivariant and all the correspondences of Fargues-Scholze we mentioned are functorial, we have the Breuil-Kisin-Fargues module produced in this way carries a semilinear $G_K$-action that commutes with all other structures of it. In this paper, we will make this correspondence explicit and call it the Fargues-Fontaine-Scholze functor $\eta_{FFS}$. We have the following result.
\begin{thm}(Theorem~\ref{FFSfunctorffandessimage})
The Fargues-Fontaine-Scholze functor
$$
\eta_{FFS}: \MF^{wa}_{K}(\varphi,N,G_K) \to \BKF(G_K)^\circ
$$
is fully faithful.
\end{thm}
Here $\BKF(G_K)^\circ$ is the isogeny category of Breuil-Kisin-Fargues $G_K$-modules(cf. Definition~\ref{defBKFGK}) and $\MF^{wa}_{K}(\varphi,N,G_K)$ is the category of weakly admissible filtered $(\varphi,N,G_K)$-modules over $K$. We will see in the isogeny category, every Breuil-Kisin-Fargues module is isomorphic to a finite free Breuil-Kisin-Fargues module(cf. Remark~\ref{isogenytoffone}), so we make the following definition:

\begin{defn}
A finite free Breuil-Kisin-Fargues $G_K$-module is called \textit{arithmetic} if, up to isogeny, it is in the essential image of the Fargues-Fontaine-Scholze functor $\eta_{FFS}$.
\end{defn}

Using Theorem~\ref{thm-intro-FFfunctor}, we can have a characterization \textit{arithmetic} Breuil-Kisin-Fargues modules. Moreover, we can also characterize the essential image of $\eta_{FFS}$ on the subcategory $\MF^{wa}_{K}(\varphi,N)$ (resp. $\MF^{wa}_{K,\varphi}$) of weakly admissible filtered $(\varphi,N)$-modules (resp. weakly admissible filtered $\varphi$-modules). Recall that for a Breuil-Kisin-Fargues module $\MMM$ it admits a de Rham realization $\MMM_{\Cp}=\MMM\otimes_{\Ainf,\theta} \O_{\Cp}$ and a crystalline realization $\MMM_{\cKur}=\MMM\otimes_{\Ainf} \cKur$ where $\cKur=W(\overline{k})[\frac{1}{p}]$. We also define $\overline{B}=\Ainf[\frac{1}{p}]/\p$ after Fargues-Fontaine, here $\p=\{ [\varpi] a \,|\, \varpi \in \m_C, \, a\in \Ainf[\frac{1}{p}] \}$.

\begin{thm}(Theorem~\ref{essentialimages})\label{mainthmintrointversion} Assume the $p$-adic monodromy theorem for $p$-adic Galois representations, then we have:
\begin{enumerate}
    \item A Breuil-Kisin-Fargues $G_K$-module $\MMM$ is arithmetic if and only if $\MMM_{\Cp}[\frac{1}{p}]$ as a $\Cp$-representation of $G_K$ is $\Cp$-admissible, i.e., it is Hodge-Tate with only 0 weight.
    \item The isogeny class of a Breuil-Kisin-Fargues $G_K$-module $\MMM$ is in the essential image of $\MF^{wa}_{K}(\varphi,N)$ if and only if $\MMM$ is arithmetic and there is a $G_K$ fixed basis inside $\MMM_{\cKur}$.
    \item The isogeny class of a Breuil-Kisin-Fargues $G_K$-module $\MMM$ is in the essential image of $\MF^{wa}_{K,\varphi}$ if and only if $(\MMM\otimes \overline{B})^{G_K}$ as a $K_0$-vector space has dimension equal to the rank of $\MMM$.
\end{enumerate}
\end{thm}

\begin{rem}\label{r:1}\quad
\begin{enumerate}
\item The $p$-adic monodromy theorem for $p$-adic Galois representations is used in the proof of (1) in the above theorem, more explicitly, we are using the fact that $\MF^{wa}_{K}(\varphi,N,G_K)$ is equivalent to the category of $\Q_p$-representations of $G_K$ that are de Rham. 
\item The terminology of being \textit{arithmetic} was first introduced in the work of Howe in \cite[\S 4]{How} using Hodge-Tate modules, we can show our definition are the same by (1).
\end{enumerate}
\end{rem}

We will be more interested in the integral version of the above theorem. 
\begin{thm}(Theorem~\ref{intversionofmaintheorem})\label{intversionintroduction}
Assume the $p$-adic monodromy theorem for $p$-adic Galois representations, and let $\Rep_{\Z_p}^{dR}(G_K)$ (resp. $\Rep_{\Z_p}^{lcr}(G_K)$, resp. $\Rep_{\Z_p}^{cris}(G_K)$) be the category of de Rham (resp. log-crystalline, resp. crystalline) representations of $G_K$ over $\Z_p$-lattices, then
\begin{enumerate}
    \item There is an equivalence of $\Rep_{\Z_p}^{dR}(G_K)$ with the category of arithmetic Breuil-Kisin-Fargues $G_K$-modules.
    \item The essential image of $\Rep_{\Z_p}^{lcr}(G_K)$ of the functor in (1) are the arithmetic Breuil-Kisin-Fargues modules such that there is a $G_K$ fixed basis inside $\MMM_{\cKur}$.
    \item The essential image of $\Rep_{\Z_p}^{cris}(G_K)$ of the functor in (1) are Breuil-Kisin-Fargues $G_K$-modules such that $(\MMM\otimes \overline{B})^{G_K}$ as a $K_0$-vector space has dimension equal to the rank of $\MMM$.
\end{enumerate}
\end{thm}

\begin{rem}\quad
\begin{enumerate}
\item Using unramified descent(cf. Lemma~\ref{lemmaunrdescent}), one can show in Theorem~\ref{mainthmintrointversion} and Theorem~\ref{intversionintroduction} one can replace $G_K$ by $I_K$ and $K_0$ by $\cKur$ in the statement.
\item From the work of \cite{BMS1}, we know there is a large class of Breuil-Kisin-Fargues $G_K$-modules comes from geometry: start with a proper smooth formal scheme $\mathfrak{X}$ over $\O_K$, and let $\overline{\mathfrak{X}}$ be its base change to $\O_{\Cp}$. Then there is a $\Ainf$-cohomology theory attached to $\overline{\mathfrak{X}}$ which is functorial in $\overline{\mathfrak{X}}$, so all the $\Ainf$-cohomology groups $H^i_{\Ainf}(\overline{\mathfrak{X}})$ carry natural semi-linear $G_K$-actions that commute with all other structures. If we take the maximal free quotients of the cohomology groups, then they are all arithmetic automatically from the \'etale-de Rham comparison theorem. So being arithmetic is the same as asking an abstract finite free Breuil-Kisin-Fargues $G_K$-module to satisfy \'etale-de Rham comparison theorem. We also have \cite[Theorem 14.1]{BMS1} shows that there is a canonical isomorphism 
\[
H^i_\mathrm{crys}(\overline{\mathfrak{X}}_{\O_{\Cp}/p}/A_\mathrm{cris})[\frac{1}{p}]\cong H^i_\mathrm{crys}(\overline{\mathfrak{X}}_{\overline{k}}/W(\overline{k}))\otimes_{W(\overline{k})}B_\mathrm{cris}^+,
\]
and we will see this is equivalent to condition $(3)$ in Theorem~\ref{intversionintroduction}.
\item From Theorem~\ref{intversionintroduction}, it is natural to ask if we can compare our theory of arithmetic BKF modules with other theories in integral $p$-adic Hodge theory. In Section~\ref{sectioncompat}, we will compare our theory with Breuil-Kisin theory (cf. \cite{KisinFcrystal}) and Liu's theory of ($\varphi, \hat{G}$)-modules theory (cf. \cite{liu-notelattice}) and a recent theory of Breuil-Kisin $G_K$-modules of Gao (cf. \cite{Gao2020breuilkisin}). In the work of Liu and Gao, they both need an input of Kisin's theory, which relies on a choice of Kummer tower $K_\infty$ over $K$, however, it was observed in \cite{liu2013compatibility} and \cite{emertongee2020moduli} that there should be some compatibility of Kisin's theory for different choices of Kummer towers over $K$ for log-crystalline representations. We will see in \S~\ref{subsectionalldescent}, our theory of arithmetic BKF modules will give a very nice explanation of such phenomena.
\item Bhatt-Scholze show in \cite{BS19} that the above $\Ainf$-cohomology theory descent to a prismatic  cohomology theory. And recently there are prismatic integral $p$-adic Hodge theory as in \cite{Bhatt-Scholze_prismaticFcrystal} and \cite{DuLiu_prismaticphihatG}. We will also explore the relation of these theories with arithmetic BKF modules in \S~\ref{subsec:prismaticarearithmetic}.
\end{enumerate}
\end{rem}

\begin{rem}
The category of Breuil-Kisin-Fargues $G_K$-modules can also be considered as the generalised representations of $G_K$ as been considered \cite{Faltingssimpson1}. And arithmetic can be thought as certain small condition. 
\end{rem}

\begin{rem}
In the work of \cite{Morrow-TsujiGeneralizedrep}, Morrow-Tsuji introduce the category $\mathrm{BKF}(\mathfrak{X},\varphi)$ of relative BKF modules over a smooth formal scheme $\mathfrak{X}$ over $\Spf(\O_K)$. And we know there is a relative version of the Fargues-Fontaine curve as been introduced in \cite{KedlayaLiu1}, we hope report some generalizations of results of this paper to the relative case in our upcoming works. 
\end{rem}

\subsection{Comparisons of different integral \texorpdfstring{$p$}{p}-adic Hodge theory}
An important property of arithmetic Breuil-Kisin-Fargues $G_K$-modules is the following rigidity result.
\begin{lem}\label{rigidintro}(Lemma~\ref{rigidityofBKF})
For any two arithmetic Breuil-Kisin-Fargues $G_K$-modules $\MMM_1$ and $\MMM_2$, if their \'etale realizations are isomorphic, then $\MMM_1\simeq\MMM_2$.
\end{lem}
In classical integral $p$-adic Hodge theory has many results that show $\Rep^\ast_{\Z_p}(G_K)$ for $\ast \in \{\cris, \lcr, \dR \}$ are equivalent to $\varphi$-modules or $(\varphi,G_K)$-modules over subrings of $\Ainf$. For example, there are theory of Wach modules (cf. \cite{Bergeronwachmod}), Kisin-Ren's theory (cf. \cite{KisinRen}), Kisin modules (cf. \cite{KisinFcrystal}), $(\varphi,\widehat{G})$-modules (cf. \cite{liu-notelattice}), and Breuil-Kisin $G_K$-modules (cf. \cite{Gao2020breuilkisin}). We will show those $\varphi$-modules or $(\varphi,G_K)$-modules after a suitable base change to $\Ainf$, give rise to arithmetic Breuil-Kisin-Fargues.

Also there are prismatic integral $p$-adic Hodge theory developed by Bhatt-Scholze in \cite{Bhatt-Scholze_prismaticFcrystal}, and their relations with $(\varphi,
\hat{G})$-modules has been studied in \cite{DuLiu_prismaticphihatG}. In particular, we have the following result.

\begin{thm}[Prismatic $F$-crystals are arithmetic at $\Ainf$, Theorem~\ref{prismaticarearithmetic}]
Let $\MM_\Prism$ be a (logarithmic) prismatic $F$-crystal over $\O_K$, then let $(\Ainf, \varphi(\ker\theta))$ be the twisted infinitesimal prism that admitting an action of $G_K$, then $\MM_\Prism((\Ainf, \varphi(\ker\theta)))$ together with the $G_K$-action is an arithmetic BKF module.
\end{thm}

The rigidity result in Lemma~\ref{rigidintro} allows us to have a universal comparison result for all those theories at $\Ainf$-level. Questions of this kind was considered in \cite{liu2013compatibility} that comparing Kisin's theory and Wach's theory. And in \cite{emertongee2020moduli}, they discussed about the compatibility of Kisin's theory for different choices of uniformizers and Kummer towers for log-crystalline representations. 

\subsection{Structure of the paper} In \S~\ref{section2}, we will review the theory of Fargues-Fontaine curve and the construction of $G_K$-equivariant modifications of vector bundles over $X_{FF}$ from weakly admissible filtered $(\varphi,N,G_K)$-modules. In \S~\ref{section3}, we will review various theories of $\varphi$-modules and their relations with vector bundles over the Fargues-Fontaine curve. Then we will review Fargues's classification theory for Breuil-Kisin-Fargues modules and the relation of isogeny classes of Breuil-Kisin-Fargues modules and admissible modifications of vector bundles. In \S~\ref{section4}, we will define the Fargues-Fontaine-Scholze functor $\eta_{FFS}$ from $\MF^{wa}_{K}(\varphi,N,G_K)$ to the isogeny category of Breuil-Kisin-Fargues $G_K$-modules and prove our main result on characterization of the essential images of $\eta_{FFS}$ on $\MF^{wa}_{K}(\varphi,N,G_K)$ and on some typical subcategories of $\MF^{wa}_{K}(\varphi,N,G_K)$. In \S~\ref{sectioncompat}, we will apply our theory to show the compatibility of many existing theories in integral $p$-adic Hodge theory. 

\subsection{Notions and conventions}\label{sub:1.4} Throughout this paper, $k$ will be a perfect field in characteristic $p$. Let $K_0=W(k)[\frac{1}{p}]$ and $\O_{K_0}=W(k)$. Let $K$ be a totally ramified finite extension of $K_0$, write $\O_K$ as the ring of integers of $K$ and let $\varpi$ be a uniformizer. Let $\cKur=W(\overline{k})[\frac{1}{p}]$ and define $I_K$ be the inertia group inside $G_K$. By a compatible system of $p^n$-th roots of $\varpi$, we mean a sequence of elements $\{\varpi_n\}_{n\geq 0}$ in $\overline{K}$ with $\varpi_0=\varpi$ and $\varpi_{n+1}^p=\varpi_{n}$ for all $n$. We normalize the  $p$-adic valuation on $K$ by $v_K(p)=1$.

Define $\Cp$ as the $p$-adic completion of $\overline{K}$, there is a unique valuation $v=v_{\Cp}$ on $\Cp$ extending the $p$-adic valuation on $K$. Let $\O_{\Cp}=\{x\in \Cp | v(x)\geq 0\}$ and let $\mathfrak{m}_{\Cp}=\{x\in\Cp | v(x)>0\}$. We will have $\O_{\Cp}/\mathfrak{m}_{\Cp}=\overline{k}$. 

Let $C=\Cp^\flat$ be the tilt of $\Cp$, then by the theory of perfectoid fields, $C$ is algebraically closed of characteristic $p$, and complete with respect to a nonarchimedean norm. Let $\O_C$ be the ring of the integers of $C$, then $\O_C=\O_{\Cp}^\flat=\varprojlim_{x\mapsto x^p}\O_{\Cp}$. $C$ is also a nonarchimedean field, we will denote the valuation by $v_C$ and normalize it by $v_C(\underline{x}) =v_{\Cp}(x_0)$ for $\underline{x}=(x_0,x_1,\ldots)\in \varprojlim_{x\mapsto x^p}\O_{\Cp}$. We will also write it by $v$ if there is no confusion. We define $\m_C=\{x\in C \, | \, v_C(x) >0\}$. Define $\Ainf=W(\O_C)$, there is a Frobenius $\varphi_{\Ainf}$ acts on $\Ainf$. $\Ainf$ is equipped with a surjection $\theta:\Ainf \to \O_{\Cp}$ satisfies $\theta([x])=x^0$ for $x=(x_0,x_1,\ldots)\in \varprojlim_{x\mapsto x^p}\O_{\Cp}$. We will have the kernel of $\theta$ is principal and let $\xi$ be a generator of $\Ker(\theta)$. We will write $\tilde\xi=\varphi(\xi)$ following the notation used in \cite{BMS1}. There is a $G_K$-action on $\Ainf$ from its action on $\O_C$, one can show $\theta$ is $G_K$-equivariant. 

\subsubsection{\textbf{Notation } $\MF_{K}(\varphi,N,G_K)$} A filtered $(\varphi, N)$-module over $K$ is a finite dimensional $K_0$-vector space $D$ equipped with two maps 
$$
\varphi, N: D\to D
$$
such that
\begin{enumerate}
\item $\varphi$ is semi-linear with respect to the Frobenius $\varphi_{K_0}$;
\item $N$ is $K_0$-linear;
\item $N\varphi = p\varphi N$.
\end{enumerate}
And a decreasing, separated and exhaustive filtration on the $K$-vector space $D_K=K\otimes_{K_0} D$.

Let $L$ be a finite Galois extension of $K$ and let $L_0=W(k_L)[\frac{1}{p}]$. A filtered $(\varphi, N, \Gal(L/K))$-module over $K$ is a filtered $(\varphi, N)$-module $D'$ over $L$ together with a semilinear action of $\Gal(L/K)$ on the $L_0$-vector space $D'$, such that:
\begin{enumerate}
\item The action is semilinear with respect to the action of $\Gal(L/K)$ on $L_0$ via $\Gal(L/K)\twoheadrightarrow\Gal(k_L/k)=\Gal(L_0/K_0)$.
\item The semilinear action of $\Gal(L/K)$ commutes with $\varphi$ and $N$.
\item The filtration on $D'\otimes_{L_0} L$ is stable under the diagonal action of $\Gal(L/K)$ on $D'\otimes_{L_0} L$, i.e., it defines a filtration on $D_K:=(D'\otimes_{L_0} L)^{G_K}$. 
\end{enumerate}

If $L^\prime$ is another finite Galois extension of $K$ containing $L$, then one can show there is a fully faithful embedding of the category of filtered $(\varphi, N, \Gal(L/K))$-modules into the category of filtered $(\varphi, N, \Gal(L^\prime/K))$-modules. One defines the category of filtered $(\varphi, N, G_K)$-modules
$$\MF_{K}(\varphi,N,G_K)$$ 
to be the limit of filtered $(\varphi, N, \Gal(L/K))$-modules over all finite Galois extensions $L$ of $K$.

Let $\cKur=W(\overline{k})[\frac{1}{p}]$, and $I_K=\Gal_{\cKur}$. For any $D \in \MF_{K}(\varphi,N,G_K)$ that is define over $L$, we let $\overline{D}=D\otimes_{L_0}\cKur$, then we have there is a continuous semilinear action of $G_K$ on $\overline{D}$ such that $I_K$ acts on $\overline{D}$ with open kernel. We will have the inverse of this is also true:
\begin{lem}\label{lemmaunrdescent}
For a finite-dimensional $\cKur$ vector space $\overline{D}$ with a continuous semilinear action of $G_K$ such that the action is trivial when restricting to an open subgroup $H$ of $I_K$, we will have 
$$
\overline{D}=(\overline{D})^{G_{L}}\otimes_{L_0}\cKur
$$
for a finite extension $L$ of $K$. Moreover, $L$ can be chosen to be a totally ramified extension of $K$, i.e., we can assume $L_0=K_0$. 
\end{lem}
\begin{proof}
We can not apply Galois descent directly in this case. We give a sketch of the proof. First since the action is continuous and $W(\overline{k})$ is DVR, so fix a lattice $\Lambda$ inside $\overline{D}$, we have $\Lambda$ is $P$ stable for an open subgroup $P$ of $G_K$. Moreover, since $G_K$ is compact, a standard trick will imply that there is $G_K$ stable lattice $\Lambda_0$ in $\overline{D}$. Let $\overline{\Lambda_0}=\Lambda_0 \mod{p}$. Let $L'$ be the corresponding totally ramified extension of $\cKur$, using Krasner’s lemma, we can always choose a uniformizer $\varpi'$ of $L'$ such that $\varpi'$ is algebraic over $K$ and defines a totally ramified extension $L$ of $K$. We have $G_L$ acts on $\overline{\Lambda_0}$ via $G_L\twoheadrightarrow G_k$, and the $\overline{k}$ vector space $\overline{\Lambda_0}$ is with the discrete topology. So we can apply Galois descent to the $G_k$ semilinear action on $\overline{\Lambda_0}$ and use the exact sequence
$$
0\to \Lambda_0 \xrightarrow{p} \Lambda_0 \to \overline{\Lambda_0} \to 0
$$
one can argue via successive approximation by lifting to show $\Lambda_0^{G_L}$ is of full rank, inverting $p$ one get 
$$
\overline{D}=(\overline{D})^{G_{L}}\otimes_{L_0}\cKur.
$$
For details of the last part of the proof, one can refer to Lemma 3.2.6 of \cite{CMInotes}.
\end{proof}

By the above lemma, for any $D \in \MF_{K}(\varphi,N,G_K)$ we can always assume the underlying $\varphi$-module is defined over $\cKur$ and equipped with a continuous semilinear $G_K$-action such that the restricted action on $I_K$ has an open kernel. We will use this fact later in this paper. 

\subsubsection{\textbf{Notation } $\MF^{wa}_{K}(\varphi,N,G_K)$} We will let $\MF^{wa}_{K}(\varphi,N,G_K)$ be the subcategory of $\MF_{K}(\varphi,N,G_K)$ consisting of weakly admissible objects. And we define its subcategory $\MF^{wa}_{K}(\varphi,N)$ (resp. $\MF^{wa}_{K,\varphi})$ which is the category of weakly admissible filtered $(\varphi,N)$-modules (resp. $\varphi$-modules).

\subsubsection{\textbf{Notation} \textbf{Linearization}} Let $R$ be a ring equipped with an endomorphism $\varphi$, for a $R$-module $M$, we write $\varphi^\ast M$ to be $M\otimes_{R,\varphi} R$. 

\subsubsection{\textbf{Conventions for semistable and log-crystalline}} In this paper, we will use the notion of \textit{log-crystalline representations} instead of \textit{semistable representations} to make a difference to the semistability of vector bundles over complete regular curves. 

\subsubsection{\textbf{Conventions for Hodge-Tate weights}}\label{conventionHTw} We will use covariant functors when relate \'etale $\varphi$-modules and Galois representation, so we will assume the cyclotomic character has Hodge-Tate weight $-1$.

\subsubsection{\textbf{Conventions for modifications of vector bundles}} We will see in our definition of modifications of vector bundles, we always mean the modification is admissible, i.e., a modifications such that one of the vector bundles is semistable of slope $0$.

\medskip
\noindent
\textbf{Acknowledgments.}
This paper is based on the author's Ph.D. thesis, I would like to thank my PhD advisor Tong Liu for his advice, supply of ideas and many useful discussions. We also thank all the committee members for their valuable comments. We thank Miaofen Chen, Hansheng Diao, Laurent Fargues, Hui Gao, Yu Min and Yupeng Wang for their interests, comments and discussions on earlier versions of this article. We learnt the notion of $p$-adic Hodge structures over $\Cp$ in terms of modification of vector bundles over the Fargues-Fontaine curve from the lectures of Fargues in the \textit{2019 Winter School on Shimura Varieties and Related Topics}, we thank Laurent Fargues for his lectures and all the organizers for the hospitality. And it should be clear that this article is greatly inspired by \cite{FF} and \cite{ScholzeWeinstein}.

\section{Fargues-Fontaine curve and the Fargues-Fontaine functor}\label{section2}

In this section, we will first review the construction and properties of the Fargues-Fontaine curve over the perfectoid field $C$, and then we will recall the Fargues-Fontaine functor from $\MF^{wa}_{K}(\varphi,N,G_K)$ to the category of $G_K$-equivariant modifications of vector bundles over the Fargues-Fontaine curve.

\subsection{The \texorpdfstring{$p$}{p}-adic fundamental curve of Fargues and Fontaine}\label{subsection11}

Recall $C$ is the tilt of $\Cp$, to define the Fargues-Fontaine curve over $C$, we first review some basic period rings of Fontaine. 

Let $\Ainf=W(\O_C)$ and recall there is a canonical surjection 
$$
\theta:\Ainf \to \O_{\Cp}$$ 
such that $\theta([x])=x_0$ for $x=(x_i)\in \varprojlim_{x\mapsto x^p}\O_{\Cp}$. Let $\xi$ be a generator of $\Ker(\theta)$. Let $\BdR^+=\varprojlim_n \Ainf[\frac{1}{p}]/(\xi^n \Ainf[\frac{1}{p}])$ and recall the topology on $\BdR^+$ is the weak topology. And $\BdR=\BdR^+[\frac{1}{\xi}]$ equipped with a $\Z$-filtration $\Fil^n \BdR= \xi^n \BdR^+$. Also recall that $\Acris$ is defined as the $p$-completed PD envelope of $\Ainf \to \O_{\Cp}$. Let $\Bcris^+=\Acris[\frac{1}{p}]$. $\Acris$ and $\Bcris^+$ can be regarded as subrings of $\BdR$, and recall the element $t=\log[\epsilon]\in \Fil^1\BdR$ is well-defined where $\epsilon=(\zeta_i)_i\in C$ is defined by a compatible system of $p^n$-th roots of units with $\zeta_0=1$. Let $\Bcris=\Bcris^+[\frac{1}{t}]$. Let $\underline{\varpi}=(\varpi_i)\in C$ defined by a compatible system of $p^n$-th roots of $\varpi$ with $\varpi_0=\varpi$ a unifromizer of $\O_K$. Then 
$$
\log[\underline{\varpi}]:=\sum_{i\geq 1}\frac{-(1-[\underline{\varpi}]/\varpi)^i}{i}
$$
is well-defined in $\Fil^1 \BdR$ and we define $\Bst^+=\Bcris^+[\log[\underline{\varpi}]]$ and $\Bst=\Bcris[\log[\underline{\varpi}]]$ with the unique $\Bcris$-derivation determined by $N(\log[\underline{\varpi}])=1$. 

\begin{rem}
Here we use the convention that $N(\log[\underline{\varpi}])=1$ as \cite[\S 10.3.2]{FF} when they define $N$ on $B_{\log}$, this is also compatible with the convention used in \cite{Kisin2adic}. We will use this fact when we compare our theory of arithmetic Breuil-Kisin-Fargues modules and Breuil-Kisin theory.
\end{rem}

\begin{rem}
When we discuss potentially log-crystalline representations, there will be an issue of change of field, we might need to consider $\log[\underline{\varpi_L}]$ where $\{\underline{\varpi_L}\}$ is defined by a compatible system of $p^n$-th roots of a uniformizer $\varpi_L$ of a finite extension $L$ of $K$. We can always choose $\varpi_L=\varpi^e \mod p$ for some $e \in \N$, and we know $\underline{\varpi_L} \in C$ only depends on the classes of $\varpi_{L,n} \mod p$, so we can always choose $\underline{\varpi_L}=\underline{\varpi}^e$. So in particular $\log[\underline{\varpi_{L}}]=e\log[\underline{\varpi}]$ and they will define the same subring $\Bst$ inside $\BdR$.
\end{rem}

For every subring $A$ of $\BdR$, define $\Fil^n A=A\cap \Fil^n \BdR$ for all $n\in\Z$. Recall that the Frobenius $\varphi$ on $\Ainf$ extensions to $\Acris$, $\Bcris^+$, $\Bcris$ and $\Bst$. We define $B_e=\Bcris^{\varphi=1}$. One has $B_e$ is actually a PID, and this is actually one of the motivations to define the Fargues-Fontaine curve. For the story behind this, one can refer to \cite{prefaceofFF}. 

For the Fargues-Fontaine curve $X_{FF}:=X_{C, \Q_p}$ (here we use the notion in \cite[Definition 6.5.1.]{FF}, where they construct Fargues-Fontaine curves $X_{F, E}$ for pairs $(F,E)$ where $F$ is any perfectoid field in characteristic $p$ and $E$ is a discrete valuation field), an abstract definition of $X_{FF}$ is that $X_{FF}$ is a scheme fits into the following Cartesian diagram 
\begin{equation}\label{diagram1}
\begin{tikzcd}
X_{FF} \arrow[dr, phantom, "\ulcorner", very near start] & \Spec(B_e) \arrow[l] \\
\Spec(\BdR^+) \arrow[u] & \Spec(\BdR) \arrow[u]\arrow[l]
\end{tikzcd}  
\end{equation}
In particular, we have $X_{FF} = \Spec(B_e)\coprod \{\infty\} $ such that $X_{FF,\infty}=\BdR^+$. Fargues-Fontaine give an explicit construction
$$
X_{FF} = \Proj \oplus_{i\geq0} (\Bcris^+)^{\varphi=p^i}.
$$
Recall $t=\log [\epsilon]$ satisfies $\varphi(t)=pt$, so it is a section of $\O(1)$. Fargues-Fontaine showed that $t$ has a unique zero $\infty\in X_{FF}$, then we summarize some of the main results in \cite{FF}:

\begin{thm}\label{FF curve}
The pointed scheme $(X_{FF}, \infty)$ fits into the diagram~\ref{diagram1}, and we have:
\begin{enumerate}
\item $X_{FF}$ is a regular noetherian scheme of Krull dimension 1, or an abstract regular curve in the sense of Fargues and Fontaine.
\item $X_e=X_{FF}\backslash \{\infty\}$ is an affine scheme $\mathrm{Spec}(B_e)$. 
\item Vector bundles $\E$ over $X_{FF}$ are equivalence to $B$-pairs $(M_e, M^+_{\mathrm{dR}},\iota)$, where $M_e=\Gamma(X_e,\E)$ is a finite projective module over $B_e$, $M^+_{\mathrm{dR}}$ is a finite free module over $\BdR^+$, and $\iota$ is an isomorphism of $M_e$ and $M^+_{\mathrm{dR}}$ over $\BdR$. And the functor is given by 
$$
\E \mapsto (\Gamma(\mathrm{Spec}(B_e),\E),\E_\infty,\iota).
$$
\end{enumerate}
\end{thm}

One advantage of having an explicit definition of $X_{FF}$ is that one can have the following construction of vector bundles from isocrystals.
\begin{thm}(Theorem 8.2.10 in \cite{FF})\label{vectorbundleandisocrystal}
Let $(D,\varphi)$ be an isocrystal over $\overline{k}$, then $(D,\varphi)$ defines a vector bundle $\E(D,\varphi)$ over $X_{FF}$ which is associated with the graded module $$ \oplus_{n\geq 0} (D\otimes_{  \cKur} \Bcris^+)^{\varphi=p^n}.$$ Moreover, this functor induces a bijection of isomorphism classes.
\end{thm}

\begin{defn}\label{definitionofslope}\quad
\begin{enumerate}
\item Let $\E$ be a vector bundle over $X_{FF}$, assume $\E\cong \E(D,\varphi)$ under the above theorem, let the multi-set $\{-\lambda_i\}$ be the slope of $(D,\varphi)$ under the Dieudonn\'e-Manin classification theorem, we define the \textit{slope} of $\E$ to be the multi-set $\{\lambda_i\}$. 
\item $\E$ is called semistable of slope $\lambda$ if and only if $\E$ corresponds a semisimple isocrystal of slope $-\lambda$. Rank $1$ vector bundle of slope $n$ is denoted by $\O(n)$ which corresponds to $(  \cKur, p^{-n}\varphi_{  \cKur})$.
\item Let $\Bun_{X_{FF}}$ be the category of vector bundles over $X_{FF}$ and let $\Bun_{X_{FF}}^{\lambda}$ be the subcategory of semistable vector bundles of slope $\lambda$.
\end{enumerate}
\end{defn}

A consequence of Theorem~\ref{vectorbundleandisocrystal} is
\begin{cor}\label{slope0}\quad
\begin{enumerate}
    \item The category of isocrystals semisimple of slope $-\lambda$ is equivalent to $\Bun_{X_{FF}}^{\lambda}$.
    \item In particular, when $\lambda=0$, the category of finite-dimensional $\Q_p$-vector spaces is equivalent to $\Bun_{X_{FF}}^{0}$, and the functor is given by 
    $$
    V \to V\otimes_{\Q_p} \O_X
    $$
    with quasi inverse
    $$
    \E \to H^0(X_{FF},\E).
    $$
\end{enumerate} 
\end{cor}
\begin{proof}
This is \cite[Theorem 9.2.2]{FF}.
\end{proof}

\begin{defn}
Let $\Modif_{X_{FF}}$ be the category of triples $(\E_0,\E_1,\iota)\in \Modif_{X_{FF}}$, where
\begin{itemize}
    \item $\E_0 \in \Bun^0_{X_{FF}}$;
    \item $\E_1 \in \Bun_{X_{FF}}$;
    \item $\iota: \E_0 |_{X_{FF}-\{\infty\}} \xrightarrow{\sim} \E_1|_{X_{FF}-\{\infty\}}$.
\end{itemize}
We will also denote this by $\E_0 \xdashrightarrow[]{} \E_1$ if no confusion arises.
\end{defn}

\begin{lem}\label{HTpairs}
$\Modif_{X_{FF}}$ is equivalent to the category of pairs $(V,\Xi)$, where 
\begin{itemize}
    \item $V$ is a finite dimensional vector space over $\Q_p$;
    \item $\Xi \in V\otimes_{\Q_p} \BdR$ is a $\BdR^+$-lattice.
\end{itemize}
\end{lem}
\begin{proof}
By (3) in Theorem~\ref{FF curve}, use the $B$-pair description of vector bundles over $X_{FF}$, we have for $\E_0 \xdashrightarrow[]{} \E_1 \in \Modif_{X_{FF}}$, then $\E_0$ and $\E_1$ have the same $B_e$-part, so the modification can can be determined by $(\E_0,\Xi)$ where $\Xi=\E_{1,\infty}$ is the $\BdR^+$ part of $\E_1$. Now by Corollary~\ref{slope0}, this is the same as $(V,\Xi)$ as in the statement of this lemma if we let $V=H^0(X_{FF},\E_0)$.
\end{proof}

\subsection{\texorpdfstring{$p$}{p}-adic representations and modifications of vector bundles on the curve}\label{subsection22}
We have $\O_X=\O_{X_{FF}}$ has a continuous action of $G_K$ via its action on $\Bcris^+$, and $\infty$ is fixed by $G_K$.
\begin{defn}\
\begin{enumerate}
    \item Let $\Bun_X(G_K)$ be the category of $G_K$-equivariant vector bundles over $X_{FF}$ such that the $G_K$-semilinear action is continuous. Here the notion of continuity can be interpreted using the $B$-pair description of vector bundles over $X_{FF}$(c.f. (3) in Theorem~\ref{FF curve}). That is continuity is the same as the $B_e$ and $\BdR^+$ parts of the vector bundle both carry a continuous semilinear action of $G_K$. Objects in $\Bun_X(G_K)$ are also called $\O_X$-representations of $G_K$.
    \item Let $\Modif_X(G_K)$ be the category of $G_K$-equivariant modifications of vector bundles, i.e. triples $(\E_0,\E_1,\iota)\in \Modif_{X_{FF}}$ such that both $\E_0,\E_1 \in \Bun_X(G_K)$ and $\iota$ is $G_K$-equivariant.
\end{enumerate}
\end{defn}

Here is another interpretation of the $G_K$-action using cocycles, cf. \cite[Proposition 9.1.5]{FF}.
\begin{lem}\label{cocyle}
Fix $\E\in\Bun_{X_{FF}}$, given a $G_K$-equivariant structure on $\E$ is the same as choosing a continuous 1-cocycles $f=(f_g)_g\in Z^1(G_K,\mathrm{Aut}(\E))$. And the isomorphic classes of $G_K$-equivariant structure on $\E$ is bijective to $H^1(G_K,\Aut(\E))$.
\end{lem}

\begin{defn}\cite[Definition 9.1.6]{FF}\label{defnwedgeaction}
Given $\E\in \Bun_X(G_K)$ and a continuous 1-cocycles $a=(a_g)_g$ inside $Z^1(G_K,\mathrm{Aut}(\E))$. And let $f=(f_g)_g\in Z^1(G_K,\mathrm{Aut}(\E))$ be a cocylce represents the $G_K$-equivariant structure on $\E$, and let $[f]\in H^1(G_K,\Aut(\E))$ be it cohomology class. Then one can show $a_g \circ [f_g] $ is a well-defined element in $H^1(G_K,\Aut(\E))$ and we will define the corresponded $\O_X$-representation by $a \wedge \E \in \Bun_X(G_K)$, i.e., it is the vector bundle with the same underlying $\O_X$-structure as $\E$ and the $G_K$-equivariant structure is defined by $a_g \circ [f_g] $ via Lemma~\ref{cocyle}.
\end{defn}

In the rest of this subsection, we review Fargues-Fontaine's functor 
$$
\eta_{FF}: \MF^{wa}_{K}(\varphi,N,G_K) \to \Modif_X(G_K).
$$ 
The main reference is \cite[\S 10.3.2]{FF}. Keep the notions as in Theorem~\ref{vectorbundleandisocrystal}.

\subsubsection{Construction of $\E(D,\varphi,N,G_K)$} Let $D$ be a filtered $(\varphi, N, G_K)$-module, recall as we have mentioned in Lemma~\ref{lemmaunrdescent}, the underlying $\varphi$-module $(D,\varphi)$ can always be treated as an isocrystal over $\overline{k}$, and let 
$$\E(D,\varphi) = \widetilde{\oplus_{n\geq 0} (D\otimes_{{  \cKur}} \Bcris^+)^{\varphi=p^n}}$$
be the vector bundle over $X_{FF}$ corresponds to $(D,\varphi)$ under the equivalence in Theorem~\ref{vectorbundleandisocrystal}. We have the $G_K$-action on $D$ defines a continuous $G_K$-action on $\E(D,\varphi)$ via the diagonal action on $D\otimes_{  \cKur} \Bcris^+$ and we use $\E(D,\varphi,G_K)$ to refer this $\O_X$-representation. 

Note that this construction is functorial, so the relation $N\varphi=p \varphi N$ tells that $N$ defines a $G_K$-equivariant map 
$$
N:\quad \E(D,\varphi,G_K) \to \E(D, p\varphi,G_K)= \E(D, \varphi,G_K)\otimes_{\O_X} \O(-1).
$$ 
Recall we define $\log[\underline{\varpi}]:=\log(\frac{[\underline{\varpi}]}{{\varpi}})$ and $t=\log[\epsilon]$ inside $\BdR^+$. For every $g\in G_K$, let $\log_{\underline{\varpi},g}=g(\log[\underline{\varpi}])-\log[\underline{\varpi}]$.  We will have $\log_{\underline{\varpi},g}=\log([\epsilon^{c_g}])=c_gt$ for a 1-cocycle $c=(c_g)_g$ valued in $\Q_p^\times$, so $\log_{\underline{\varpi},g}$ is inside $(\Bcris^+)^{\varphi=p}$, i.e. $\log_{\underline{\varpi},g}$ defines a morphism from $\O(-1)$ to $\O$. So we know the composition:
\[
\begin{tikzcd}
\beta_g: \E(D,\varphi) \arrow[r, "N"] & \E(D, \varphi)\otimes \O(-1) \arrow[r, "\mathrm{Id}\otimes \log_{\underline{\varpi},g}"] & \E(D,\varphi)\otimes_{\O_X} \O_X=\E(D,\varphi)
\end{tikzcd}
\]
defines an element in $\mathrm{End}(\E(D,\varphi))$. Moreover, since $\log_{\underline{\varpi},g}$ is a $\Q_p$ multiple of $t$, so $(\log_{\underline{\varpi},g})_g$ defines an element in $Z^1(G_K, (\Bcris^+)^{\varphi=p})$, so $\beta=(\beta_g)_g$ defined as above is an element in $Z^1(G_K, \mathrm{End}(\E(D,\varphi)))$, and by the nilpotence of $N$, we have the image of $\beta$ lies in the nilpotent elements of $\mathrm{End}(\E(D,\varphi))$. So we can define
\begin{equation}\label{alphaaction}
\alpha=(\alpha_{g})_{g}=(-\exp(\beta_g))_g \in Z^1(G_K, \mathrm{Aut}(\E(D,\varphi))).
\end{equation}
Fargues-Fontaine define a new $G_K$-action on the vector bundle $\E(D,\varphi)$ by twisting, we let
$$
\E(D,\varphi,N,G_K)=\alpha \wedge \E(D,\varphi,G_K).
$$
i.e., $\E(D,\varphi,N,G_K)$ is isomorphic to $\E(D,\varphi)$ as vector bundle, and the $G_K$-action is given by twisting the $G_K$-action of $\E(D,\varphi,G_K)$ with the 1-cocycle $\alpha$(cf. Definition~\ref{defnwedgeaction}). 

There is another way of thinking about this new action.
\begin{lem}\label{Dhat}
Let $D$ be a $(\varphi,N)$-module. There is an $\varphi$-equivariant isomorphic
\begin{IEEEeqnarray*}{+rCl+x*}
D\otimes_{{  \cKur}} \Bcris^+ & \quad \xrightarrow{\sim} \quad & (D\otimes_{{  \cKur}} \Bst^+)^{N=0} \\
 y \quad & \mapsto & \hat{y}:=\sum_{i\geq 0} \frac{(-1)^i}{i!} N^i(v) \otimes \log([\underline{\varpi}])^i,
\end{IEEEeqnarray*}
for every $y\in D$ and extends linearly to $D\otimes_{{  \cKur}} \Bcris^+$.
\end{lem}
\begin{proof}
Here we will use a similar computation in \cite[\S 7.2]{liu-notelattice}. Write $X=\log([\underline{\varpi}])$ and $\gamma_i(X)=\frac{X^i}{i!}$. We have $N(\gamma_i(X))=\gamma_{i-1}(X)$ by our convention $N(X)=1$. And we have
$$
N(\hat{y})=\sum_{i\geq 0} {(-1)^i} (N^{i+1}(v) \otimes \gamma_i(X) + N^{i}(v) \otimes \gamma_{i-1}(X)) =0
$$
i.e., $\hat{y}\in (D\otimes_{{  \cKur}} \Bst^+)^{N=0}$. A direct computation also shows that $y \to \hat{y}$ is $\varphi$-equivariant. And we define
$$
\widehat{D}=\{\hat{v}\,| \,v\in D\}.
$$
One can check $\widehat{D}$ has a structure of isocrystal over $\overline{k}$ via 
$$\cKur\xhookrightarrow{} (\Bst^+)^{N=0} \xhookrightarrow{} \Bst^+
$$ 
and it is isomorphic to $D$. Moreover via $\Bcris^+ = (\Bst^+)^{N=0} \xhookrightarrow{} \Bst^+$, one have $(D\otimes_{{  \cKur}} \Bst^+)^{N=0}$ is a $\Bcris^+$-module which is isomorphic to $\widehat{D}\otimes_{\cKur} \Bcris^+$.  
\end{proof}

If $D\in \MF_{K}(\varphi,N,G_K)$, given $ D\otimes \Bst^+$ the diagonal action of $G_K$, then we can view $D\otimes \Bcris^+$ as a $G_K$-stable subspace of $ D\otimes \Bst^+$ via the the above lemma.
\begin{lem}\label{hataction}
With the notions in Lemma \ref{Dhat}, we have $G_K$ acts on $(D\otimes_{{  \cKur}} \Bst^+)^{N=0}$ by
$$
g(\hat{y})= \exp(-\log_{\underline{\varpi},g}N))(\widehat{g(y)}).
$$
\end{lem}
\begin{proof}
We can write $\hat{y}=\exp(-\log([\underline{\varpi}])N)(y)$ for $y\in D$. Use the fact that the $G_K$-action commutes with $N$, we have 
\begin{IEEEeqnarray*}{+rCl+x*}
g(\hat{y}) & = & g(\exp(-\log([\underline{\varpi}])N)(y)) \\
           & = & \exp(-g(\log([\underline{\varpi}]))N)(g(y)) \\
           & = & \exp((-\log_{\underline{\varpi},g}-\log([\underline{\varpi}]))N)(g(y)) \\
           & = & \exp(-\log_{\underline{\varpi},g}N)\exp((-\log([\underline{\varpi}])N)(g(y)) \\
           & = & \exp(-\log_{\underline{\varpi},g}N)(\widehat{g(y)}).
\end{IEEEeqnarray*}
\end{proof}

In terms of vector bundles, the map $y \mapsto \hat{y}$ in Lemma \ref{Dhat} defines an isomorphism of graded modules:
$$
\hat{h} :\oplus_{i\geq 0} (D\otimes \Bcris^+)^{\varphi=p^i} \xrightarrow{\sim} \oplus_{i\geq 0} (D\otimes \Bst^+)^{N=0,\varphi=p^i},
$$
So it defines an automorphism $\hat{h}_\E$ of $\E=\E(D)$. Given $\oplus_{i\geq 0} (D\otimes \Bst^+)^{N=0,\varphi=p^i}$ the $G_K$-action defined in Lemma \ref{hataction}, then the Lemma implies this action, in terms of cocylces, is given by $(\alpha_g \circ \hat{h}_\E \circ f_g \circ \hat{h}_\E^{-1})_g$, where $(f_g)_g$ is a coclyce defines the $G_K$-structure of $\E(D,\varphi,G_K)$ and $\alpha_g$ is defined in \ref{alphaaction}. In other words, this action is exactly $\alpha\wedge \E(D,\varphi,G_K)$ by Definition \ref{defnwedgeaction}.

\subsubsection{$B_e$ and $\BdR^+$-representations of $G_K$}
From (3) in Theorem~\ref{FF curve}, one can also describe an $\O_X$-representation in term of $B$-pairs, we want to note that the terminology of $B$-pairs was first appeared in the work of Berger \cite{BergerBpairs}.
\begin{prop} \label{Bpairs}
An $\O_X$-representation $\E$ is equivalence to a $G_K$-representation on a $B$-pair $(M_e, M^+_{\mathrm{dR}},\iota)$. Here $M_e=\Gamma(X_e,\E)$ (resp. $M^+_{\mathrm{dR}}$) is a finite free $B_e$-module (resp. $\BdR^+$-module) with a continuous $G_K$-semilinear action, and $\iota$ is a $G_K$-equivariant isomorphism of $M_e$ and $M^+_{\mathrm{dR}}$ over $\BdR$.
\end{prop}

We want to describe $\E(D,\varphi,N,G_K)$ in terms of $B$-pairs. For the $\BdR^+$ part, we want to mention the following results. They play a critical role in our theory. 

\begin{defn}
A $\BdR$(resp. $\BdR^+$)-representation of $G_K$ is a finite free $\BdR$ (resp. $\BdR^+$) module together with a continuous $G_K$-semilinear action. We say a $\BdR$ (resp. $\BdR^+$)-representation $W$ of $G_K$ is flat (resp. generically flat) if $\dim_K(W^{G_K})=\dim_{\BdR}W$ (resp. $\dim_K(W[\frac{1}{t}]^{G_K})=\dim_{\BdR}W[\frac{1}{t}]$).
\end{defn}

\begin{prop}\label{BdRflat}\
\begin{enumerate}
    \item $-\otimes_K \BdR$ induces an equivalence of finite dimensional $K$ vector space and flat $\BdR$-representations of $G_K$, and the quasi-inverse is given by $W \mapsto W^{G_K}$.
    \item $(V, \Fil^\bullet V) \mapsto \Fil^0(V\otimes_K \BdR)$ induces an equivalence of filtered $K$ vector spaces and generically flat $\BdR^+$-representations of $G_K$, and the quasi-inverse is given by $W \mapsto (W[\frac{1}{t}]^{G_K},(t^\bullet W)^{G_K})$.
    \item A rank $d$ $\BdR^+$-representation $W$ of $G_K$ is generically flat if and only if there exist $(a_1,\ldots, a_d)\in \Z^d$ such that 
    $$
    W \simeq \oplus_{i=1}^d t^{a_i}\BdR^+.
    $$
\end{enumerate}
\end{prop}

\begin{proof}
This is \S 10.4 \cite{FF}. 
\end{proof}

In this paper, we will also consider $\BdR^+$-flat representations.
\begin{defn}
A $\BdR^+$-representation $W$ of $G_K$ of rank $d$ is flat if there is a $G_K$-equivariant isomorphism 
    $$
    W \simeq \oplus_{i=1}^d \BdR^+.
    $$
\end{defn}

\begin{prop}\label{BdR+flat}\
\begin{enumerate}
    \item The $\BdR^+$-representation $W$ of $G_K$ is $\BdR^+$-flat if and only if $W\otimes_{\BdR^+} \Cp$ as a $\Cp$-representation is $\Cp$-admissible, i.e., it is Hodge-Tate with all weight equal to $0$.
    \item There is an equivalence of $\BdR^+$-flat representations of $G_K$ and finite dimensional $K$-vector space, and the functor is given by $W \mapsto (W[\frac{1}{t}])^{G_K}$ with quasi-inverse given by $V \mapsto V\otimes_{K}\BdR^+$.
\end{enumerate}
\end{prop}

\begin{proof}
For (1), the author's proof first appeared in \cite[Proposition F.13.]{emertongee2020moduli}. Assume that $W \otimes_{\BdR^+} \Cp$ has a $G_K$-stable basis, we need to construct a $G_K$-stable basis of $W$. Recall that $\BdR^+$ has the weak topology such that the quotients $\BdR^+/(\xi^n)\simeq \Cp^n$ are equipped with the $p$-adic topology, so we will argue via constructing a successive lifting of $G_K$-fixed basis $\{e^{(n)}_i\}$ of $W/t^i W$ satisfying $e^{(n)}_i=e^{(n+1)}_i \mod (\xi^n)$ for all $n\geq 1$.

For $n=1$, this is given by the assumption. And assume we have already constructed $\{e^{(n-1)}_i\}$, we choose an arbitrary lifting $\{\tilde{e}^{(n)}_i\}$ of $\{e^{(n-1)}_i\}$ to $W\otimes_{\BdR^+} \BdR^+/(\xi^n)$. For $g\in G_K$, assume $g$ acts by 

$$g.(\tilde{e}^{(n)}_1,\tilde{e}^{(n)}_2,...,\tilde{e}^{(n)}_d)=(\tilde{e}^{(n)}_1,\tilde{e}^{(n)}_2,...,\tilde{e}^{(n)}_d). A_g^{(n)}$$

we have the following diagram:
\[
\begin{tikzcd}
g.(\tilde{e}^{(n)}_1,\tilde{e}^{(n)}_2,...,\tilde{e}^{(n)}_d) \arrow[r, equal] \arrow[d, twoheadrightarrow, "\mod (\xi^{n-1})" ] & (\tilde{e}^{(n)}_1,\tilde{e}^{(n)}_2,...,\tilde{e}^{(n)}_d)A_g^{(n)} \arrow[d, twoheadrightarrow, "\mod (\xi^{n-1})" ] \\
g.(e^{(n-1)}_1, e^{(n-1)}_2,..., e^{(n-1)}_d) \arrow[r, equal ] & (e^{(n-1)}_1, e^{(n-1)}_2,..., e^{(n-1)}_d)
\end{tikzcd}
\]

Let $d=\rank(\M)$, we have $A_g^{(n)}$ is a $d\times d$ matrix with coefficients in $\BdR^+/(\xi^n)$ satisfies $A_g^{(n)}\equiv I_d \mod(\xi^{n-1})$. Recall, $t$ is also a generator of $(\xi)$ in $\BdR^+$, so one can write $$A_g^{(n)}=I_d + t^{n-1} B_g^{(n)}$$ for some $B_g^{(n)}\in M_{d}(\BdR^+/(\xi^{n}))$ and this expression is uniquely determined by the class of $B_g^{(n)}$ in $M_{d}(\BdR^+/(\xi))$.

The map $g \mapsto A_g^{(n)}$ defines a 1-cocycle of $G_K$ in $\GL_d(\BdR^+/(\xi^{n}))$. And by a simple computation, we will have $B_{hg}^{(n)}=B_h^{(n)}+\chi^{n-1}(h)h(B_g^{(n)})$, i.e., the map $g \mapsto B_g^{(n)}$ defines a 1-cocycle of $G_K$ in $M_{d}(\BdR^+/(\xi)(n-1))$, that is, a 1-cocycle of $G_K$ in $M_{d\times d}(\Cp(n-1))$. Since $n>1$ by our assumption, by Tate-Sen, we have such a 1-cocycle is a 1-coboundary, i.e., there is a $\alpha^{(n)}\in M_{d}(\Cp)$ such that $B_g^{(n)}=\chi^{n-1}(g)g(\alpha^{(n)})-\alpha^{(n)}$.

Then let $(e^{(n)}_1, e^{(n)}_2,..., e^{(n)}_d) = (\tilde{e}^{(n)}_1,\tilde{e}^{(n)}_2,...,\tilde{e}^{(n)})(I_d - t^{n-1} \alpha^{(n)})$, we have: 
\begin{IEEEeqnarray*}{+rCl+x*}
g.(e^{(n)}_1,..., e^{(n)}_d) & = & (\tilde{e}^{(n)}_1,...,\tilde{e}^{(n)}_d)(I_d+ t^{n-1}B_g^{(n)})g.(I_d-t^{n-1}\alpha^{(n)}) \\
    & = & (\tilde{e}^{(n)}_1,...,\tilde{e}^{(n)}_d)(I_d+ t^{n-1}B_g^{(n)}-g.t^{n-1}g.(\alpha^{(n)})) \\
    & = & (\tilde{e}^{(n)}_1,...,\tilde{e}^{(n)}_d)(I_d+ t^{n-1}(\chi^{n-1}(g)g(\alpha^{(n)})-\alpha^{(n)})-\chi^{n-1}(g)t^{n-1}g.(\alpha^{(n)})) \\
    & = & (\tilde{e}^{(n)}_1,...,\tilde{e}^{(n)}_d)(I_d - t^{n-1} \alpha^{(n)})
\end{IEEEeqnarray*}
and satisfies $e^{(n-1)}_i \equiv e^{(n)}_i \mod (\xi^{n-1})$.

For (2), it is a consequence of (2) in Proposition~\ref{BdRflat}.
\end{proof}

\begin{rem}\label{remarkonEG}\
\begin{enumerate}
    \item Given a $p$-adic de Rham representations $V$ of $G_K$, there are two ways to assign a $G_K$-stable $\BdR^+$-lattice in $V\otimes_{\Q_p}\BdR$. Using (2) of Proposition~\ref{BdRflat}, we know there is a lattice $\Xi_0$ that corresponds to $D_{\dR}(V)$ with the Hodge filtration. On the other hand, we can also give $D_{\dR}(V)$ the trivial filtration, and produce another $G_K$-stable lattice $\Xi_1=\D_\dR(V)\otimes_{K} \BdR^+$.
    \item $\BdR^+$ is not ($\Q_p,G_K$)-regular since $t\BdR^+$ is $G_K$-stable but $t$ is not a unit. However one can define for a $\Q_p$-representation $V$ it is called $\BdR^+$-admissible if the following holds
    $$
(V\otimes_{\Q_p} \BdR^+)^{G_K} \otimes_K \BdR^+ \xrightarrow{\sim} V \otimes_K \BdR^+.
$$
Then one can show by Proposition~\ref{BdRflat}, $\BdR^+$-admissibility is equivalent to $\BdR$-admissibility plus the condition that the Hodge-Tate weights are non-negative (one can also refer to \cite[Exercise 15.5.5]{CMInotes}). In \cite[Proposition F.13.]{emertongee2020moduli}, the proof actually shows that the representation $T$ is $\BdR^+$-admissibility, and that is the reason why they require the Frobenius on Breuil-Kisin-Fargues modules is an endomorphism. We will come back to this issue.
\end{enumerate} 
\end{rem}

Recall the following result of Berger on $B_e$-representations and $(\varphi,N)$-modules.
\begin{prop}\label{Bereps}\cite[10.3.20]{FF}
There are functors 
\begin{IEEEeqnarray*}{+rCl+x*}
\mathscr{D}_{\log} :  \Rep_{B_e}(G_K) & \xrightarrow{} \quad & (\varphi,N)\text{-modules} \\
  M & \mapsto & (M\otimes_{B_e} \Bst)^{G_{K'}}
\end{IEEEeqnarray*}
and
\begin{IEEEeqnarray*}{+rCl+x*}
\mathscr{V}_{\log} :  (\varphi,N)\text{-modules} & \xrightarrow{} \quad &  \Rep_{B_e}(G_K)\\
 D \quad & \mapsto & (D\otimes_{\cKur} \Bst)^{\varphi=1,N=0}
\end{IEEEeqnarray*}
We have $\mathrm{Id}\xrightarrow{\sim} \mathscr{D}_{\log} \circ \mathscr{V}_{\log}$, and $\mathscr{V}_{\log} \circ \mathscr{D}_{\log}(M)=M$ if and only if $M$ is in the essential image of $\mathscr{V}_{\log}$. Similar result holds when replacing $(\varphi,N)$-modules by $\varphi$-modules and $\Bst$ by $\Bcris$.
\end{prop}

\subsubsection{$B$-pairs for $\E(D,\varphi,N,G_K)$ and $\E(D,\varphi,N,\Fil^\bullet,G_K)$}Fargues-Fontaine give a nice description of the $\O_X$-representation $\E(D,\varphi,N,G_K)$ in terms of $B$-pair.

\begin{prop}\label{Bpairs of E1}
Let $(M_e, M^+_{\mathrm{dR}},\iota)$ be the $G_K$-representation on a $B$-pair corresponds to the $\O_X$-representation $\E(D,\varphi,N,G_K)$, then
$$
M_e \xrightarrow{\sim} \mathscr{V}_{\log}(D)=(D\otimes_{\cKur} \Bst)^{\varphi=1,N=0}
$$
and 
$$
M^+_{\mathrm{dR}}\xrightarrow{\sim} D_K\otimes_K \BdR^+
$$
as $G_K$-modules. Here $D\otimes_{\cKur} \Bst$ is given by the diagonal action of $G_K$ and $D_K=(D\otimes_{\cKur} \BdR)^{G_K}$.
\end{prop}

\begin{proof}
For the $B_e$ part, it follows from \cite[\S 10.3.5]{FF}. For the $\BdR^+$ part, one can follow the argument of \cite[Proposition 10.3.18]{FF}. Or we can use Proposition~\ref{BdR+flat}, it is equivalence to show it is $\BdR^+$-flat. But the latter can be easily seen from the fact 
$$
\exp(-\log_{\underline{\varpi},g}N)\equiv 1 \mod{\Ker(\theta)}
$$
since $\log_{\underline{\varpi},g}=c_g t$ is a multiply of $t$. 
\end{proof}

Now let us construct a modification of $\E(D,\varphi,N,\Fi^\bullet,G_K)$. By Proposition~\ref{BdRflat} a filtration on $D_K=(D\otimes_{\cKur} \BdR)^{G_K}$ is equivalent to a $G_K$ stable $\BdR^+$ lattice in $D\otimes_{\cKur} \BdR$. Fargues-Fontaine define the $\O_X$-representation $\E(D,\varphi,N,\Fil^\bullet,G_K)$ by letting:
$$
\E(D,\varphi,N,\Fi^\bullet,G_K)|_{X_{FF}\backslash\infty}=\E(D,\varphi,N,G_K)|_{X_{FF}\backslash\infty}
$$
and 
$$
\E(D,\varphi,N,\Fi^\bullet,G_K)_\infty=\Fi^0(D_K\otimes_{K} \BdR).
$$
To show $\E(D,\varphi,N,\Fil^\bullet,G_K) \xdashrightarrow[]{} \E(D,\varphi,N,G_K)$ defines a $G_K$-equivariant modification, it remains to show the admissibility, i.e., $\E(D,\varphi,N,\Fi^\bullet,G_K)$ is semistable of slope $0$. This follows from the Proposition below, which is one of the key result of Fargues-Fontaine explaining the relation of weakly admissibility and the slope of vector bundles.

\begin{prop}\label{Bpairs of E0}
The filtered $(\varphi, N, G_K)$-module $D$ is weakly admissible if and only if $\E(D,\varphi,N,\Fi^\bullet,G_K)$ is semistable of slope 0. Moreover, let $V$ be the potentially log-crystalline representation of $G_K$ corresponding to $D$, then there is a $G_K$-equivariant isomorphism $$V=H^0(X_{FF},\E(D,\varphi,N,\Fi^\bullet,G_K)),$$ 
where $V$ is the potentially log-crystalline representation corresponding to the data $(D,\varphi,N,\Fi^\bullet,G_K)$.
\end{prop}

\begin{proof}
For the first part, it is stated in \cite[\S 10.5.3, Remark 10.5.8]{FF} for filtered $\varphi$-modules, and for filtered $(\varphi,N)$-modules, it is \cite[Proposition 5.6]{prefaceofFF}. For the second part of the proof, let $V$ be the potentially log-crystalline representation corresponding to $(D,\varphi,N,\Fi^\bullet,G_K)$ and let $\E_V=V\otimes_{\Q_p} \O_X$ be the corresponding slope 0 $\O_X$-representation of $G_K$. It is equivalent to show
$$\E_V=\E(D,\varphi,N,\Fi^\bullet,G_K)$$ 
and we can prove it by comparing the $B$-pairs of $\E(D,\varphi,N,\Fi^\bullet,G_K)$ and $\E_V$ by of Proposition~\ref{Bpairs}. While we have the $B_e$-part of the $\O_X$-representation $\E(D,\varphi,N,\Fi^\bullet,G_K)$ is the same as the $B_e$-part of the $\O_X$-representation $\E(D,\varphi,N,G_K)$ by construction, and which is equal to 
$$\mathscr{V}_{\log}(D)=(D\otimes B_{st})^{\varphi=1,\,N=0}$$
by Proposition~\ref{Bpairs of E1}. The $\BdR^+$-part of $\E(D,\varphi,N,\Fi^\bullet,G_K)$ is the $\BdR^+$-representation
$$
\Fi^0(D_{K}\otimes_{K} \BdR)
$$
by definition. 

On the other hand, the $B$-pair correspond to $\E_V$ is
$$
(V\otimes_{\Q_p}B_e,V\otimes_{\Q_p}\BdR^+).
$$
Since $V$ is potentially log-crystalline, so we have $V\otimes_{\Q_p}B_e$ is potentially log-crystalline as a $B_e$-representation in the sense that there is a $G_K$-equivariant isomorphism
$$
V\otimes_{\Q_p}B_e=\Big(\big((V\otimes_{\Q_p}B_e)\otimes_{B_e}B_{\st}\big)^{G_L}\otimes_{L_0} B_{\st}\Big)^{\varphi=1,\,N=0}
$$
for a finite Galois extension $L$ of $K$ by Proposition 10.3.20 of \cite{FF}, and $((V\otimes_{\Q_p}B_e)\otimes_{B_e}B_{\st}\big)^{G_L}\simeq D$ as $(\varphi,N,G_K)$-modules. For the $\BdR^+$-part, since $V$ is de Rham, so the $\BdR^+$-representation $V\otimes_{\Q_p}\BdR^+$ is \textit{generically flat}, so Proposition~\ref{BdRflat} shows that there is a $G_K$-equivariant isomorphism
$$
V\otimes_{\Q_p}\BdR^+ = \Fi^0(D_{\mathrm{dR}}(V)\otimes_{K} \BdR).
$$
\end{proof}

\subsubsection{The Fargues-Fontaine functor}
\begin{defn}
We define 
\begin{IEEEeqnarray*}{+rCl+x*}
\eta_{FF} :  \MF_{K}^{wa}(\varphi,N,G_K) & \quad \xrightarrow{} \quad & \Modif_X(G_K) \\
 D \quad & \mapsto & \E(D,\varphi,N,\Fil^\bullet,G_K) \xdashrightarrow[]{} \E(D,\varphi,N,G_K).
\end{IEEEeqnarray*}
\end{defn}

Let's show $\eta_{FF}$ is fully faithful.
\begin{lem}\label{lemmaGKTHpairs}
$\Modif_X(G_K)$ is equivalent to the category of pairs $(V,\Xi)$ where $V$ is a representation of $G_K$ over a $\Q_p$-vector space and $\Xi$ is a $G_K$ stable $\BdR^+$-lattice in $V\otimes_{\Q_p} \BdR^+$.
\end{lem}

\begin{proof}
Directly from Lemma~\ref{HTpairs}.
\end{proof}

\begin{thm}\label{etaFFisff}
$\eta_{FF}$ is fully faithful, and the essential image of $\eta_{FF}$ in terms of $(V,\Xi)$ such that $V$ is a de Rham representation of $G_K$ and $\Xi$ is a $\BdR^+$ stable lattice in $V\otimes \BdR$ such that as $\BdR^+$-representations of $G_K$, it is flat.
\end{thm}
\begin{proof}
Using Lemma~\ref{lemmaGKTHpairs} and Proposition~\ref{Bpairs}, we can write down $\eta_{FF}$ in terms of $(V,\Xi)$-pairs:
\begin{IEEEeqnarray*}{+rCl+x*}
\tilde{\eta}_{FF} :  \MF_{K}^{wa}(\varphi,N,G_K) & \quad \xrightarrow{} \quad & \{(V,\Xi)\} \\
 D \quad & \mapsto & (V(D), D_{\mathrm{dR}}(V(D))\otimes_{K} \BdR^+)
\end{IEEEeqnarray*}
where $V(D)$ is the potentially log-crystalline representation of $G_K$ corresponds to $D$. So it is obvious that $\eta_{FF}$ is fully faithful and the essential image lies inside the category defined in the theorem. 

And if a pair $(V,\Xi)$ such that $V$ is a de Rham representation of $G_K$ and $\Xi$ a flat $\BdR^+$-representations of $G_K$ inside $V\otimes \BdR$, then by the equivalence in Proposition~\ref{BdRflat}, we will have $\Xi=D_{\mathrm{dR}}(V)\otimes_{K} \BdR^+$, by the $p$-adic monodromy theorem of $p$-adic Galois representations, we have $V$ is potentially log-crystalline, so corresponds to some $D\in \MF_{K}(\varphi,N,G_K)$. And again by the computations in Proposition~\ref{Bpairs of E0}, we have the the $B$-pair corresponds to $\eta_{FF}(D)$ is $(V,\Xi)$.
\end{proof}

\section{\texorpdfstring{$\varphi$}{f}-modules and Breuil-Kisin-Fargues modules}\label{section3}
This section will discuss $\varphi$-modules over various base rings and their relation with vector bundles over $X_{FF}$. And we will also discuss Breuil-Kisin-Fargues modules and recall their relations with modifications of vector bundles over the Fargues-Fontaine curve.

\subsection{\texorpdfstring{$\varphi$}{f}-modules and vector bundles} Recall we have the rings in characteristic $p$: $k$, $\overline{k}$, $\O_C$, and $C$ on which the Frobenius endomorphism is defined. Moreover, since we have
$$
\O_C=\varprojlim_{x\to x^p} \O_{\overline{K}}/p = \varprojlim_{x\to x^p} \O_{\overline{K}_0}/p.
$$
So there is a canonical map $k \to \O_C$ comes from $k\simeq \O_{K_0}/p \to \O_{\overline{K}_0}/p$ which has a uniquely lifting to a section $\overline{k} \to \O_C$ of $\O_C \to \overline{k}$ by \'etaleness. This defines a canonical $\varphi$-equivariant section 
$$
W(\overline{k}) \to \Ainf.
$$

Let's recall the following constructions of rings with $\varphi$-endomorphisms, the main reference of this part is \cite[\S 11]{FF}. Define $B^+=\cap_n \varphi^n(\Bcris^+)$ and let $B$ be the completion of $B^b:=\Ainf[\frac{1}{[\underline{\varpi}]},\frac{1}{p}]$ with respect to $\{\lvert \cdot \rvert_\rho\}_{\rho\in (0,1)}$, here for $\rho\in (0,1)$, define
$$
\lvert \sum_{i \gg -\infty} [x_i]p^i  \rvert_\rho = \sup\{ \lvert x_i \rvert \rho^i \} 
$$
where $\sum_{i \gg -\infty} [x_i]p^i \in B^b$. For $\rho\in(0,1)$, let $B^+_\rho$ be the completion of $\Ainf[\frac{1}{p}]$ with respect to $\lvert \cdot \rvert_\rho$, then $B^+=\cap_{\rho\in (0,1)} B^+_\rho$ and $B^+ \xhookrightarrow{} B$. Consider the ideal $\p=\cup_n(\varphi^{-n}[\underline{\varpi}]) \subset \Ainf[\frac{1}{p}]$ where $\varpi$ is any pseudo uniformizer in $C$, and define $\tilde{\p}=\p B^+$, we will have $\Ainf[\frac{1}{p}] \to B^+$ induces an isomorphism $\Ainf[\frac{1}{p}]/\p = B^+/\tilde{\p}$ by \cite[11.1.1]{FF}, and let $\overline{B}$ denote this quotient. Then we have $\overline{B}$ is a local domain with maximal ideal $\m_{\overline{B}}:=W(\m_C)[\frac{1}{p}]/\p$ and residue field $\cKur$. Let's recall the following properties of $B^+$:
\begin{exam}\label{Brho}
For $\rho\in (0,1)$, we have $B^+_{\rho}=\widehat{\Ainf[\frac{[a]}{p}]}[\frac{1}{p}]$, where $\lvert a \rvert=\rho$ and the completion is for the $(p,[\varpi])$-adic topology for any $\varpi\in \m_C$, but we have $[a]$ is divisible by $p$ in $\widehat{\Ainf[\frac{[a]}{p}]}$, so the topology is the same as the $p$-adic topology. In particular, we have $\widehat{\Ainf[\frac{[a]}{p}]}$ is also completed under the $[\varpi]$-adic topology.
\end{exam}

\begin{defn}
Let $R$ be a ring that is equipped with a Frobenius endomorphism $\varphi$, a $\varphi$-module over $R$ is a finite projective module $M$ together with 
$$
\varphi^\ast M:= M\otimes_{R,\varphi} R \xrightarrow{\sim} M.
$$
And we use $\phimod_R$ to denote the category of $\varphi$-modules over $R$. For a $\varphi$-module $M$, we define its global section $H^0(M)$ by
$$
H^0(M)=M^{\varphi=1}.
$$
\end{defn}

\begin{prop}\label{strucutesofphimod}\
\begin{enumerate}
    \item $M \mapsto M \otimes_{B^+} \overline{B} $ induces an equivalence of $\phimod_{B^+}$ and $\phimod_{\overline{B}}$.
    \item For $(d,h) \in \N\times \N_{>0}$ with $(h,d)=1$, and $R\in\{\cKur,B,B^+,\Bcris^+,\overline{B}\}$, let $R(\frac{d}{h})$ be the $\varphi$-module free on the bases $\{e_i\}_{i=0}^{h-1}$ and $\varphi(e_i)=e_{i+1}$ for $i\in [0,h-2]$ and $\varphi(e_{h-1})=p^d e_{0}$. Then any $M\in \phimod_{R}$ is a direct sum of $R(\frac{d}{h})$ over pairs $(d,h)$.
    \item By (2), for any $M\in \phimod_{B^+}$ (resp. $M_{\overline{B}}\in \phimod_{\overline{B}}$), let $\breve{M}=M \otimes_{B^+} \breve{K}$ (resp. $\breve{M}=M_{\overline{B}} \otimes_{\overline{B}} \breve{K}$), then there is a $\varphi$-equivariant section $s: \breve{M} \to M$ (resp. $s: \breve{M} \to M_{\overline{B}}$) reducing to the identity over $\cKur$.
\end{enumerate}
\end{prop}

\begin{proof}
This is \cite[Theorem 11.1.7]{FF}. 
\end{proof}

\begin{lem}\label{Bdh}\
\begin{enumerate}
    \item For all $h\in \N_{>0}$ and $d\in \Z$, the natural maps between $B$, $B^+$, $\Bcris^+$ and $\overline{B}$ induces
$$
B^{\varphi^h=p^d} =(B^+)^{\varphi^h=p^d} = (\Bcris^+)^{\varphi^h=p^d}=\overline{B}^{\varphi^h=p^d}.
$$
\item If $d<0$, then we will have
$$
B^{\varphi^h=p^d} = 0.
$$
\end{enumerate}
\end{lem}
\begin{proof}
(2) is from \cite[Proposition 4.1.2]{FF}. For (1), in $loc.cit.$, they show $B^{\varphi^h=p^d} =(B^+)^{\varphi^h=p^d}$ and $(B^+)^{\varphi^h=p^d}=\overline{B}^{\varphi^h=p^d}$ when $d>0$. We give a different approach to prove this for all $d$.

Let $R(\frac{d}{h})$ be as in (2) in in Proposition~\ref{strucutesofphimod} and for $R\in\{\overline{B},B^+,B,\Bcris^+\}$, we have 
$$
R^{\varphi^h=p^d}=\Hom_{\varphi}(R(\frac{d}{h}),R).
$$
Since the categories of $\varphi$-modules over $R$ are all equivalent to each other by \cite[\S 11.1]{FF}. So the above $\Hom$ are all equal with each other.
\end{proof}

\begin{prop}\label{propostionMdh}
Let $M^+\in \phimod_{B^+}$ and define $M=M^+\otimes_{B^+} B$ (resp. $M_\cris=M^+\otimes_{B^+} \Bcris^+$, resp. $M_{\overline{B}}=M^+\otimes_{B^+} \overline{B}$). Then we have for all $h\in \N_{>0}$ and $d\in \Z$.
$$
(M^+)^{\varphi^h=p^d} =M^{\varphi^h=p^d} = (M_\cris)^{\varphi^h=p^d}=M_{\overline{B}}^{\varphi^h=p^d}.
$$
\end{prop}

\begin{proof}
Use the same trick as in Lemma~\ref{Bdh}, we have $M^{\varphi^h=p^d}=\Hom_{\varphi}(B(\frac{d}{h}),M)$, and similar formula for $M^+$ $M_\cris$, and $M_{\overline{B}}$. So the fully faithfulness of the base change functors implies the result.
\end{proof}

\begin{thm}\cite{FF}\label{t:5}
Let the category $\phimod_{B^+}$, $\phimod_{B}$, $\phimod_{\Bcris^+}$, and $\phimod_{\overline{B}}$ are all equivalent to the category of $\Bun_{X_{FF}}$. And one direction of the functor is
$$
(M,\varphi)\mapsto \E(M) = \widetilde{\oplus_{n\geq 0} {M}^{\varphi=p^n}}
$$
for $(M,\varphi) \in \phimod_{B^+}$ (resp. $\phimod_{B}$, resp. $\phimod_{\Bcris^+}$, resp. $\phimod_{\overline{B}}$).
\end{thm}

\begin{proof}
The equivalence of $\phimod_{B^+}$, $\phimod_{B}$ and $\Bun_{X_{FF}}$ is from \cite[\S 11.4]{FF}. And $M \mapsto M\otimes_{B^+}\Bcris^+$ induces an equivalence of $\phimod_{B^+}$ and $\phimod_{\Bcris^+}$ and with quasi-inverse functor given by $ M \mapsto \cap_n \varphi^n(M)$. The equivalence of $\phimod_{B^+}$ and $\phimod_{\overline{B}}$ was mentioned in \cite[\S 11]{FF}.

For the last statement for the consistency of the functor from different $\varphi$-modules to $\Bun_{X_{FF}}$, first we have by \cite[Theorem 11.1.9]{FF}. The functor 
$$
(M,\varphi)\mapsto \widetilde{\oplus_{n\geq 0} {M}^{\varphi=p^n}}
$$
induces an equivalence of $\phimod_{B^+}$ and $\Bun_{X_{FF}}$. We want to emphasize that for $M^+$, $M$, $M_\cris$, and $M_{\overline{B}}$ are as in Proposition~\ref{propostionMdh}, we have 
$$
(M^+)^{\varphi=p^n} \simeq M^{\varphi=p^n} \simeq (M_\cris)^{\varphi=p^n} \simeq M_{\overline{B}}^{\varphi=p^n}
$$
for all $n$, so define same graded modules. And this is precisely what we have proved in Proposition~\ref{propostionMdh}.
\end{proof}

\begin{defn}
With the notions in Proposition~\ref{strucutesofphimod}, and let $\lambda=\frac{d}{h}$. We define the category $\phimod_{R}^{\lambda}$ of semisimple $\varphi$-modules of slope $\lambda$ over $R$ to be the subcategory of $\phimod_{R}$ consisting $M$ that is finite a direct sum of $R(\lambda)$.
\end{defn}

From Theorem~\ref{t:5} and Corollary~\ref{slope0}, we have the following lemma:
\begin{lem}\label{canonicalsection}
For $R\in\{\cKur,B,B^+,\Bcris^+,\overline{B}\}$, $\phimod_{R}^{\lambda}$ are canonically equivalent to each other under the natural base changes. And they are all equivalent to the category of vector bundles semistable of slope $-\lambda$ over $X_{FF}$.  

In particular, let $M \in \phimod_{R}^{\lambda}$ with $R\in\{B,B^+,\Bcris^+,\overline{B}\}$ and let $M_{\cKur}$ be its base change to $\cKur$, then there is a canonical $\varphi$-equivariant section $s: M_{\cKur} \to M$ that reduces to identity over $\cKur$.
\end{lem}

The ring $\overline{B}$ plays an important role in our theory, so we give a brief discussion of it here.

\begin{defn}
Let $B^+_{(1,1]}=\cup_{\rho\in(0,1)}B^+_{\rho}$, and define $\p_1=\p B^+_{(1,1]}$ and $\m_1=W(\m_C)B^+_{(1,1]}$. It is easy to see from Example~\ref{Brho} that $B^+_{(1,1]}/\p_1=\overline{B}$ and $B^+_{(1,1]}/\m_1=\cKur$. This $B^+_{(1,1]}$ carries a Frobenius endomorphism induces from the Frobenius endomorphisms on $B^+_\rho$, so we can define $\phimod_{B^+_{(1,1]}}$. We give $B^+_{(1,1]}$ the final topology induced by the inclusions $B^+_{\rho} \to B^+_{(1,1]}$.
\end{defn}

\begin{lem}\label{B11andB+}
The bases change functor along the natural map $B^+ \to B^+_{(1,1]}$ induces an equivalence of $\phimod_{B^+}$ and $\phimod_{B^+_{(1,1]}}$.
\end{lem}
\begin{proof}
The proof is similar to \cite[Proposition 12.3.5]{ScholzeWeinstein}. Let $(M,\varphi_M)\in \phimod_{B^+_{(1,1]}}$, then there is $\rho\in(0,1)$, such that $(M,\varphi_M)$ is defined over $B^+_\rho$. Then we can use the fact $\varphi: B^+_{\rho} \xrightarrow{\sim} B^+_{\rho^p}$ to extend this $\varphi$-module over $B^+_\rho$ to a $\varphi$-module over $B^+$ using $\varphi_M$. 
\end{proof}

\begin{lem}
We have $(B^+_{(1,1]},\p_1)$ is a henselian pair in the sense of \cite[\href{https://stacks.math.columbia.edu/tag/09XI}{Tag 09XI}]{stacks-project}. 
\end{lem}
\begin{proof}
Here we are using Gabber's definition of henselian pair following \cite[\href{https://stacks.math.columbia.edu/tag/09XI}{Tag 09XI}]{stacks-project}, we will show:\\
$\p_1$ is contained in the Jacobson radical of $B^+_{(1,1]}$ and every monic polynomial $f(T)\in B^+_{(1,1]}[T]$ of the form
$$
f(T)=T^n(T-1)+a_nT^n+\ldots+a_1T+a_0
$$
with $a_n,\ldots,a_0\in\p_1$ and $n\geq 1$ has a root in $1+\p_1$.

Let $x\in \p_1$ and $b\in B^+_{(1,1]}$, then $xb=[\varpi]\frac{y}{p^n}$ for a pseudo uniformizer $\varpi \in \m_C$ and $y\in \widehat{\Ainf[\frac{[a]}{p}]}$ from the discussion in Example~\ref{Bdh}. Then we have $\frac{y}{p^n}$ is inside $\widehat{\Ainf[\frac{[a]}{p^{n+1}}]}=\widehat{\Ainf[\frac{[\tilde{a}]}{p}]}$. Since $\widehat{\Ainf[\frac{[\tilde{a}]}{p}]}$ is $[\varpi]$-adically complete, so $1+xb$ is a unit, in particular, $x$ is in the Jacobson radical of $B^+_{(1,1]}$.

Similarly, for the monic polynomial
$$
f(T)=T^n(T-1)+a_nT^n+\ldots+a_1T+a_0
$$
we can find $\varpi_0,a_0 \in \m_C$ such that $a_i\in [\varpi_0]\widehat{\Ainf[\frac{[a_0]}{p}]}$ for all $i$. Again, using the fact that $\widehat{\Ainf[\frac{[a_0]}{p}]}$ is $[\varpi_0]$-adically complete. We have the pair $(\widehat{\Ainf[\frac{[a_0]}{p}]},[\varpi_0])$ is henselian by \cite[\href{https://stacks.math.columbia.edu/tag/0ALJ}{Tag 0ALJ}]{stacks-project}. So $f(T)$ has a root in $1+[\varpi_0]\widehat{\Ainf[\frac{[a_0]}{p}]}\subset 1+\p_1$.
\end{proof}

\begin{cor}
The reduction from $B^+_{(1,1]}$ to $\overline{B}$ induces an equivalence of $\phimod_{B^+_{(1,1]}}$ and $\phimod_{\overline{B}}$.
\end{cor}
\begin{proof}
This is from of \cite[Lemma 4.1.26]{ALB}. First we want to point out the prove of Lemma 4.1.26 in $loc.cit.$ does not use the prism structure in their statement. To apply their lemma, we just need to show $\varphi$ is topologically nilpotent on $\p_1$. For $x\in\p_1$, $x=[\varpi]\frac{a}{p^n}$ for some $\varpi \in \m_C$ and $a\in B^+_{\rho}$ satisfies $\lvert a \rvert_\rho \leq 1$, so 
$$
\lvert\varphi^n([\varpi]\frac{a}{p^n})\rvert_\rho=\frac{1}{p^n}\lvert[\varpi]^{p^n}\rvert_\rho \lvert\varphi^n(a)\rvert_\rho
$$
converge to $0$.
\end{proof}

\begin{rem}
The above argument gives an different proof of (1) in Proposition~\ref{strucutesofphimod} by Lemma~\ref{B11andB+}. 
\end{rem}

\begin{rem}
Let's also gives a geometric interpretation of the ring $\overline{B}$. Recall we have the following picture 

\tikzset{every picture/.style={line width=0.75pt}} 
\[
\begin{tikzpicture}[x=0.75pt,y=0.75pt,yscale=-1,xscale=1]

\draw    (162.09,223.68) -- (162.09,45.8) ;
\draw [shift={(162.09,43.8)}, rotate = 450] [color={rgb, 255:red, 0; green, 0; blue, 0 }  ][line width=0.75]    (10.93,-3.29) .. controls (6.95,-1.4) and (3.31,-0.3) .. (0,0) .. controls (3.31,0.3) and (6.95,1.4) .. (10.93,3.29)   ;
\draw    (162.09,223.68) -- (355.9,223.68) ;
\draw [shift={(357.9,223.68)}, rotate = 180] [color={rgb, 255:red, 0; green, 0; blue, 0 }  ][line width=0.75]    (10.93,-3.29) .. controls (6.95,-1.4) and (3.31,-0.3) .. (0,0) .. controls (3.31,0.3) and (6.95,1.4) .. (10.93,3.29)   ;
\draw    (162.09,223.68) -- (289.01,117.04) ;
\draw    (239.48,138.38) .. controls (230.39,116.68) and (229.28,117.9) .. (207.62,106.77) ;
\draw [shift={(205.91,105.89)}, rotate = 387.43] [color={rgb, 255:red, 0; green, 0; blue, 0 }  ][line width=0.75]    (10.93,-3.29) .. controls (6.95,-1.4) and (3.31,-0.3) .. (0,0) .. controls (3.31,0.3) and (6.95,1.4) .. (10.93,3.29)   ;
\draw   (156.03,223.68) .. controls (156.03,220.35) and (158.74,217.65) .. (162.09,217.65) .. controls (165.44,217.65) and (168.15,220.35) .. (168.15,223.68) .. controls (168.15,227.01) and (165.44,229.71) .. (162.09,229.71) .. controls (158.74,229.71) and (156.03,227.01) .. (156.03,223.68) -- cycle ;
\draw    (163.01,20.04) .. controls (230.88,20.59) and (290.53,40.02) .. (330.7,80.35) .. controls (363.33,113.1) and (383.11,159.64) .. (384.01,221.04) ;
\draw    (301.01,106.04) -- (320.52,88.38) ;
\draw [shift={(322.01,87.04)}, rotate = 497.86] [color={rgb, 255:red, 0; green, 0; blue, 0 }  ][line width=0.75]    (10.93,-3.29) .. controls (6.95,-1.4) and (3.31,-0.3) .. (0,0) .. controls (3.31,0.3) and (6.95,1.4) .. (10.93,3.29)   ;
\draw    (330.7,80.35) -- (334.66,77.04) ;
\draw    (384.01,221.04) -- (389.66,221.04) ;
\draw    (163.01,20.04) -- (162.66,14.04) ;

\draw (247.48,157.84) node [anchor=north west][inner sep=0.75pt]   [align=left] {$x_{\Cp}$};
\draw (243.52,225.8) node [anchor=north west][inner sep=0.75pt]   [align=left] {$x_{\text{\'et}}$};
\draw (90,116.39) node [anchor=north west][inner sep=0.75pt]   [align=left] {$x_{\cris}=x_{\overline{B}}$};
\draw (134.6,223.02) node [anchor=north west][inner sep=0.75pt]   [align=left] {$x_{\overline{k}}$};
\draw (201.84,121.96) node [anchor=north west][inner sep=0.75pt]   [align=left] {$\varphi$};
\draw (359.9,226.68) node [anchor=north west][inner sep=0.75pt]   [align=left] {$p=0$};
\draw (140.76,29.24) node [anchor=north west][inner sep=0.75pt]   [align=left] {$[\underline{p}]=0$};
\draw (300,10) node [anchor=north west][inner sep=0.75pt]   [align=left] {$\mathcal{Y}=\Spa(\Ainf)\backslash \{x_{\overline{k}}\}$};
\draw (288,76) node [anchor=north west][inner sep=0.75pt]   [align=left] {$\kappa$};
\draw (164,-1) node [anchor=north west][inner sep=0.75pt]   [align=left] {$\infty$};
\draw (338,59) node [anchor=north west][inner sep=0.75pt]   [align=left] {$1$};
\draw (395,211) node [anchor=north west][inner sep=0.75pt]   [align=left] {$0$};
\end{tikzpicture}
\]
Here $\mathcal{Y}$ is defined by the locus in $\Spa(\Ainf)$ such that $p(x)\neq 0$ or $[\underline{p}](x)\neq 0$, there is a continuous map $\kappa: \mathcal{Y} \to [0,\infty]$ given by the relative position of $x$ to $[\underline{p}]$ and $p$. And for $x\in \mathcal{Y}$, $\kappa(x)=\infty$ implies $x$ factor through $\Ainf \to \Ainf[\frac{1}{p}]\to \Ainf[\frac{1}{p}]/\p=\overline{B}$. Here $\Ainf \to \Ainf[\frac{1}{p}]$ comes from the condition $p(x)\neq 0$ (since we require $[\underline{p}](x)=0$) and $\p$ is the ideal subject to the condition $\{[\underline{p}](x)=0\}$. In particular, $x_{\cris}$ in this picture should be considered not just a point but the space $\overline{B}$. The space $\overline{B}$ is actually very huge in the sense that it has infinite Krull dimension (cf. \cite{LLAinf}). Moreover, one can show there is an uncountable chain of prime ideals inside $\overline{B}$ that is fixed by $\varphi$ (cf. \cite{HDAinf}).

We will see $\varphi$-modules over $\overline{B}$ plays an important role when studying the crystalline Galois representations. 
\end{rem}

\subsection{Breuil-Kisin-Fargues modules and modifications of vector bundles}
The main reference of this section is \cite{BMS1} and \cite{ansch_BKFCM}. 

\begin{defn}\
\begin{enumerate}
    \item A Breuil-Kisin-Fargues module is a finite presented module $\MMM$ over $\Ainf$ with an isomorphism
$$\varphi_{\MMM}\,:\,\MMM\otimes_{\Ainf,\varphi}\Ainf[\frac{1}{\tilde{\xi}}]\simeq \MMM[\frac{1}{\tilde{\xi}}]$$
such that $\MMM[\frac{1}{p}]$ is a finite free $\Ainf[\frac{1}{p}]$-module. Recall here $\tilde{\xi}=\varphi(\xi)$ as we defined in \S~\ref{subsection11}. 
    \item A Breuil-Kisin-Fargues module is a finite free if the underlying module is finite free over $\Ainf$. In this paper, we define $\BKF$ to be the category of \textbf{finite free} Breuil-Kisin-Fargues modules. And we let $\BKF^\circ$ to be the isogeny category of $\BKF$. We will also use BKF modules to refer objects in $\BKF$. 
    \item A Breuil-Kisin-Fargues module $\MMM$ is called effective if $\varphi_{\MMM}(\MMM\otimes_{\Ainf,\varphi}\Ainf) \subset \MMM$.
\end{enumerate}
\end{defn}

\begin{exam}\label{effective}
There are free rank $1$ Breuil-Kisin-Fargues modules $\Ainf\{n\}$ defined for all $n \in \Z$(cf. \cite[Example 4.24]{BMS1}). Just note in this paper, our $\Ainf\{1\}$ is the $\Ainf\{-1\}$ in \cite[Example 4.24]{BMS1} due to our conventions on Hodge-Tate weights~\ref{conventionHTw}, by this convention we will have $\Ainf\{n\}\otimes_{\Ainf} \overline{B} = \overline{B}(n)$ defined in (2) Proposition~\ref{strucutesofphimod}. For any $\MMM\in \BKF$, we define 
$$
\MMM\{n\}=\MMM\otimes_{\Ainf} \Ainf\{n\}.
$$
It can be shown that $\MMM\{n\}$ is effective for $n\gg 0$(cf. \cite[\S 2.3.6.]{CorIris}). 
\end{exam}

\begin{rem}\label{isogenytoffone}
In the proof of \cite[Lemma 3.9]{ansch_BKFCM}, any Breuil-Kisin-Fargues module is isogeny to a finite free one, so $\BKF^\circ$ is also the isogeny category of \textbf{all} Breuil-Kisin-Fargues modules. 
\end{rem}

\begin{defn}
Let $\HT_{\Z_p}$ be the category of pairs $(T,\Xi)$, where $T$ is a finite free $\Z_p$-lattice and $\Xi$ is a $\BdR^+$-lattice in $T\otimes_{\Z_p}\BdR$. And let $\HT$ be the isogeny category of $\HT_{\Z_p}$, i.e., $\HT$ consists of pairs $(V,\Xi)$ where $V$ is a finite $\Q_p$-vector space and $\Xi$ is a $\BdR^+$-lattice in $V\otimes_{\Q_p}\BdR$.
\end{defn}

\begin{defn}
Given $\MMM\in \BKF$, we define:
\begin{enumerate}
    \item $T(\MMM)=(\MMM\otimes_{\Ainf}W(C))^{\varphi=1}$;
    \item $V(\MMM)=(\MMM\otimes_{\Ainf}W(C))^{\varphi=1}[\frac{1}{p}]$;
    \item $\MMM_{\Cp}=\MMM\otimes_{\Ainf,\theta} {\Cp}$;
    \item $\MMM_{\cris}=\MMM\otimes_{\Ainf} \Bcris^+$;
    \item $\MMM_{\breve{K}}=\MMM\otimes_{\Ainf} W(\overline{k})[\frac{1}{p}]$.
\end{enumerate}
\end{defn}

\begin{rem}
It is easy to see that $V(\MMM)$, $\MMM_{\Cp}$, $\MMM_{\cris}$ and $\MMM_{\breve{K}}$ are well-defined on the isogeny classes of BKF modules.
\end{rem}

\begin{lem}
For $\MMM\in \BKF$, we have $(T(\MMM), \MMM\otimes_{\Ainf}\BdR^+) \in \HT_{\Z_p}$.
\end{lem}
\begin{proof}
By \cite[Lemma 4.26]{BMS1}, $T(\MMM)$ is a $\Z_p$-lattice of the same rank equal to $\MMM$, and $T(\MMM)\otimes_{\Z_p} \BdR = \MMM\otimes \BdR$, so $(T(\MMM), \MMM\otimes_{\Ainf}\BdR^+) \in \HT$. 
\end{proof}

\begin{lem}
For $\MMM\in \BKF$, we have $\MMM_{\cris}\in \phimod_{\Bcris^+}$
\end{lem}
\begin{proof}
It is enough to show $\tilde{\xi}$ is a unit in $\Bcris^+$. There are many way to see this, for readers familiar with the language of $\delta$-rings \cite{BS19}, we give a short proof as follow, we have $\Acris=\Ainf\{\frac{\tilde{\xi}}{p}\}^\wedge_{\delta}$, and both $\tilde{\xi}$ and ${p}$ are distinguished, so $\frac{\tilde{\xi}}{p}$ is a unit in $\Acris$.
\end{proof}

\begin{thm}(Fargues' classification theorem, cf. \cite[Theorem 4.28]{BMS1})\label{Farguesclassification}
The functor $$\MMM \mapsto (T(\MMM), \MMM\otimes_{\Ainf}\BdR^+)$$ induces an equivalence of $\BKF$ and $\HT_{\Z_p}$.
\end{thm}

\begin{cor}
There is an equivalence of $\Modif_{X_{FF}}$ and $\BKF^\circ$.
\end{cor}

\begin{proof}
By Lemma~\ref{HTpairs} and Theorem~\ref{Farguesclassification}, both category are equivalent to $\HT$. 
\end{proof}

\begin{rem}
One can regard a modification $\E_0 \xdashrightarrow[]{} \E_1$ as a modification of $\E_0$ by a $\BdR^+$-lattice at a formal neighborhood at $\infty$, and this gives the equivalence in Lemma~\ref{HTpairs}. On the other hand, one can also view it as a modification of $\E_1$ by a lattice in $\E_{1,\infty}\otimes_{\BdR^+}\BdR$. And this gives the following result. 
\end{rem}

\begin{prop}\label{BcrisofBKF}
The functor
$$
\MMM \to (\MMM_\cris, V(\MMM)\otimes_{\Q_p} \BdR^+).
$$
defines a fully faithfully functor from $\BKF^\circ$ to $\{(M_{\Bcris^+},\Xi)\}$, where the latter is the category consisting of pairs $(M_{\Bcris^+},\Xi)$, with $M_{\Bcris^+}\in \phimod_{\Bcris^+}$ and $\Xi$ is a $\BdR^+$ lattice in $M\otimes_{\Bcris^+} \BdR$.
\end{prop}

\begin{proof}
If $\MMM$ corresponds to $\E_0\xdashrightarrow[]{} \E_1$, then if one carefully tracks the functor defined by Scholze in \cite[Theorem 14.1.1]{ScholzeWeinstein}, one has $\E_1$ corresponds to the $\varphi$-module $\MMM_{\cris}$. And we have $\E_0$ is a semistable slope $0$ vector bundle with global section equals to $V(\MMM)$, in particular, the $\BdR^+$-lattice defined by $\E_{0,\infty}$ is $V(\MMM)\otimes_{\Q_p} \BdR^+$. By the $B$-pair description of vector bundles over $X_{FF}$ we have the functor is fully faithful.
\end{proof}

We also want to discuss about sections to the natural reduction of $\MMM$ to $\MMM\otimes_{\Ainf} {\cKur}$.
\begin{lem}\label{lemmasection}
For $\MMM\in \BKF$, let $M_{\cris}=\MMM\otimes \Bcris^+$, and $M_{\overline{B}}=\MMM\otimes_{\Ainf} \overline{B}$. For any perfect subfield $l \subset \overline{k}$ and let $L_0=W(l)[\frac{1}{p}]$. For any $\varphi$-stable $L_0$-vector space $V\subset M_{\overline{B}}$ such that $V \otimes_{L_0} \overline{B} = M_{\overline{B}}$, there is an unique $\varphi$-equivariant section $s: V\to M_{\Bcris^+}$ such that modulo $\p\Bcris^+$ we get identity on $V$. 
\end{lem}
\begin{proof}
The proof is combination of \cite[Lemma 4.5.6]{Cais-Liu} and \cite[Proposition 4.26]{Farguesconjecture}. Since we will use the formula of the section $s$, let us give an explicit construction of the section following the idea in $loc.cit.$.

First, we can always reduce to the case that $\MMM$ is effective, since we can always twist $\MMM$ by $\Ainf\{n\}$ and $V$ by $L_0(n)$ simultaneously by some $n\gg 0$ to ensure that $\MMM$ and $V$ being effective by Example~\ref{effective}. So we can assume there is a $\varphi$-stable $\O_{L_0}$-lattice $\Lambda$ inside $V$. Since we have $\MMM_{\cKur}=(\Lambda \otimes_{\O_{L_0}}\overline{B})\otimes_{\overline{B}} \cKur$, we can view $\Lambda$ as a sub-$\varphi$-module of $\MMM[\frac{1}{p}] \otimes_{\Ainf[\frac{1}{p}]} \overline{B}$.

Let $\{\overline{e}_i\}_{i=1}^d$ be a basis of $\Lambda$ and choose an arbitrary lifting $\{e_1,\ldots, e_d\}$ as a basis of $\MMM\otimes_{\Ainf} \Ainf[\frac{1}{p}]$, we can always replace $\{\overline{e}_i\}_{i=1}^d$ by $\{p^k\overline{e}_i\}_{i=1}^d$ and $\{e_1,\ldots, e_d\}$ by $\{p^k e_1,\ldots, p^k e_d\}$ to ensure $\{e_1,\ldots, e_d\}$ is inside $\MMM$. We define $A_0 \in M_d(O_{L_0})\subset \overline{B}$ be the matrix defined by
$$
\varphi(\overline{e}_1,\ldots,\overline{e}_d)=(\overline{e}_1,\ldots,\overline{e}_d)A_0.
$$
let $A \in M_d(\Ainf)$ be the matrix defined by
$$
\varphi({e}_1,\ldots,{e}_d)=({e}_1,\ldots,{e}_d)A.
$$
Then to find $s$ it is enough to find $Y\in M_d(\Bcris^+)$ such that $Y\equiv I \mod{\p}$ and
$$
Y A_0 =A \varphi(Y),
$$
and the section is defined by 
$$
s(\overline{e}_1,\ldots,\overline{e}_d)=({e}_1,\ldots,{e}_d)Y.
$$
For $n\in \N_{>0}$, define 
$$
Y_n=A\varphi(A)\cdots \varphi^{n-1}(A)\varphi^{n-1}(A_0^{-1})\cdots \varphi(A_0^{-1})A_0^{-1}.
$$
It is enough to show $Y_n$ converges to $Y\in M_d(\Bcris^+)$. Since we have $A\equiv A_0 \mod{\p}$ and there is $a\in \N$ such that $A_0^{-1}\in \frac{1}{p^a}M_d(O_{L_0})$, we have there is a $B_0\in M_d(\O_{L_0})$ such that $AB_0={p^a}I + [\varpi]Z$, where $\varpi$ is a pseudo uniformizer in $C$ and $Z\in M_d(\Ainf)$. We will have
$$
Y_n-Y_{n-1}=A\varphi(A)\cdots \varphi^{n-2}(A) \frac{[\varpi]^{p^{(n-1)}}}{{p^{an-a}}} \varphi^{n-1}(Z)\varphi^{n-2}(B_0)\cdots \varphi(B_0)B_0.
$$
It is easy to show that $\frac{[\varpi]^{p^{(n-1)}}}{{p^{an-a}}}$ converge to $0$ in $\Bcris^+$ (actually in $B^+$). So we have $Y_n$ converges to some $Y$. The uniqueness is also the same as in \cite[Lemma 4.5.6]{Cais-Liu} and we skip it here. 
\end{proof}

Apply the above lemma to the case $l=\overline{k}$, i.e., $L_0=\cKur$, we have the following Corollary on the section $s$ define in (3) of Proposition~\ref{strucutesofphimod}.
\begin{cor}\label{corosection}
For $\MMM\in \BKF$, let $M_{\cris}=\MMM\otimes \Bcris^+$, $M_{\overline{B}}=\MMM\otimes_{\Ainf} \overline{B}$, and $\MMM_{\cKur}=\MMM \otimes_{\Ainf} \cKur$. For any section 
$$
\overline{s}: \MMM_{\cKur} \to M_{\overline{B}}
$$ 
reducing to the identity over $\cKur$. There is a unique section
$$
{s}: M_{\cKur} \to M_{\cris}
$$ 
such that reducing $\p\Bcris^+$ is identify to $\overline{s}$.
\end{cor}
\begin{proof}
Just apply Lemma~\ref{lemmasection} to $\overline{s}(M_{\cKur})$.
\end{proof}

There is a definition related to the above results.

\begin{defn}(Rigidifications c.f. \cite{ansch_BKFCM})\label{rigidification}
For any $\MMM\in \BKF$ define $\MMM_{\cris}=\MMM\otimes_{\Ainf} \Bcris^+$, $\MMM_{\overline{B}}=\MMM\otimes_{\Ainf} \overline{B}$, and $\MMM_{\cKur}=\MMM\otimes_{\Ainf}\cKur$. A rigidification of $\MMM$ over $\Bcris^+$ (resp. $\overline{B}$) is a $\varphi$-equivariant section 
$$
{s}: \MMM_{\cKur} \to M_{\cris} \quad (\text{resp. }\, \overline{s}: \MMM_{\cKur} \to M_{\overline{B}})
$$ 
reducing to the identity over $\cKur$. 
\end{defn}

By Corollary~\ref{corosection}, we have.
\begin{cor}
For any $\MMM \in \BKF$, a rigidification of $\MMM$ over $\overline{B}$ is the same as a rigidification over $\Bcris^+$.
\end{cor}

\section{Fargues-Fontaine-Scholze functor and arithmetic Breuil-Kisin-Fargues modules}\label{section4}

\begin{defn}\label{defBKFGK}
Let $\BKF(G_K)$ be the category of finite free Breuil-Kisin-Fargues $G_K$-modules consists of Breuil-Kisin-Fargues module equipped with a continuous semilinear $G_K$-action that commutes with $\varphi_{\MMM}$. And let $\BKF(G_K)^\circ$ be the corresponded isogeny category. 
\end{defn}

One easily deduce the following ``$G_K$-version" of the Fargues' classification theorem.

\begin{prop}\label{GKFarguesclass}
There is an equivalence of $\Modif_X(G_K)$ and $\BKF(G_K)^\circ$.
\end{prop}

\begin{defn}
We define the Fargues-Fontaine-Scholze functor 
$$
\eta_{FFS}: \MF^{wa}_{K}(\varphi,N,G_K) \to \BKF(G_K)^\circ
$$
to be the composition of the Fargues-Fontaine functor $\eta_{FF}$ with the functor in Proposition~\ref{GKFarguesclass}.
\end{defn}

In this section, we will show the essential image of $\eta_{FFS}$ is characterized by what we call arithmatic Breuil-Kisin-Fargues $G_K$-modules; moreover, we will also characterize the essential images corresponds to the subcategory $\MF^{wa}_{K}(\varphi,N)$ (resp. $\MF^{wa}_{K,\varphi})$.

\subsection{Arithmetic Breuil-Kisin-Fargues \texorpdfstring{$G_K$}{GK} modules}
\begin{defn}
Let $\MMM\in\BKF(G_K)$, $\MMM$ is called arithmetic if and only if $\MMM_{\Cp}$ as a $\Cp$-representation of $G_K$ is $\Cp$-admissible, i.e., it is Hodge-Tate with only 0 weight. We let $\BKF^a(G_K)$ to be the category of arithmetic Breuil-Kisin-Fargues modules, and $\BKF^a(G_K)^\circ$ be the corresponded isogeny category.
\end{defn}

\begin{rem}
The condition of being arithmetic is well defined on the isogeny class of $\MMM$.
\end{rem}

\begin{thm}\label{FFSfunctorffandessimage}
The Fargues-Fontaine-Scholze functor is fully faithful. The essential image of $\eta_{FFS}$ consists of arithmetic Breuil-Kisin-Fargues modules.
\end{thm}

\begin{proof}
The fully faithfulness comes from the fully faithfullness of $\eta_{FF}$ and the equivalence in Proposition~\ref{GKFarguesclass}. By Theorem~\ref{etaFFisff}, we know the essential image of $\eta_{FF}$ consist of Hodge-Tate modules $(V,\Xi)$ such that $V$ is a de Rham representation and $\Xi=D_{\dR}(V)\otimes_{\Q_p}\BdR^+$. So it is enough to show that if $\MMM$ is an arithmetic Breuil-Kisin-Fargues $G_K$-modules, then $V=V(\MMM)$ is de Rham and 
$$\MMM\otimes \BdR^+\xrightarrow{\sim} D_{\dR}(V)\otimes_{\Q_p}\BdR^+.$$ 
First, we have $\MMM_{\Cp}$ being $\Cp$-admissible, is equivalent to $\MMM \otimes_{\Ainf} \BdR^+$ is $\BdR^+$-flat by Proposition~\ref{BdR+flat}. In particular, $\MMM \otimes_{\Ainf} \BdR^+$ is generically flat, that is $V\otimes_{\Q_p} \BdR \simeq \MMM\otimes_{\Ainf} \BdR$ is $\BdR$-flat, in other word, $V$ is de Rham. By (2) in Proposition~\ref{BdR+flat}, $\MMM \otimes_{\Ainf} \BdR^+$ corresponds to $D_\dR(V)=(V\otimes \BdR)^{G_K}$ with the trivial filtration, so we have
$$
\MMM \otimes_{\Ainf} \BdR^+ \simeq D_{\dR}(V) \otimes_{K} \BdR^+.
$$
\end{proof}

\begin{rem}\
\begin{enumerate}
    \item In \cite{emertongee2020moduli}, the definition of BKF modules is slightly different from our definition. They use $\xi$ instead of $\tilde{\xi}$ in the definition. So the our definition is differed by a Frobenius twist. Moreover, they require the $\varphi_\MMM$ is an endomorphism on $\MMM$ this will result in the Hodge-Tate weights of $T(\MMM)$ being non-negative. 
    \item In \cite[Propostion F.13]{emertongee2020moduli}, they actually show that $T$ is $\BdR^+$-admissible in the sense of (2) in Remark~\ref{remarkonEG}, which requires the non-negativity of Hodge-Tate weights of $T(\MMM)$.
\end{enumerate}
\end{rem}

\begin{rem}
As we have mentioned in Remark~\ref{remarkonEG}, for a de Rham representation $T$, a natural way to define a $G_K$-stable $\BdR^+$-lattice in $D_{\dR}(T)\otimes_K \BdR$ is given by
$$
\Xi=\Fil^0(D_{\dR}(T)\otimes_K \BdR).
$$
Here if we give $D_{\dR}(T)$ the Hodge filtration, we get a lattice $\Xi_0$ and if we give $D_{\dR}(T)$ the trivial filtration, we let $\Xi_1$ the the resulting lattice. One can check the pair $(T,\Xi_0)$ corresponds to the trivial BKF $G_K$-module $T\otimes_{\Z_p} \Ainf$ and $(T,\Xi_1)$ corresponds to the arithmetic Breuil-Kisin-Fargues modules in Theorem~\ref{FFSfunctorffandessimage}.
\end{rem}

\subsection{Conditions for log-crystalline and crystalline representations.} Let $\MF^{wa}_{K}(\varphi,N)$ (resp. $\MF^{wa}_{K,\varphi}$) be the subcategory of $\MF^{wa}_{K}(\varphi,N,G_K)$ consists of weakly admissible filtered $(\varphi,N)$-modules (resp. weakly admissible filtered $\varphi$-modules). In this section, we will characterize of the essential image of $\eta_{FFS}$ of $\MF^{wa}_{K}(\varphi,N)$ and $\MF^{wa}_{K,\varphi}$.

\begin{thm}\label{essentialimages}\
\begin{enumerate}
    \item[(1)]  An arithmetic Breuil-Kisin-Fargues $G_K$-module $\MMM$ is in the essential image of $\MF^{wa}_{K}(\varphi,N)$ if and only if there is a $G_K$-fixed basis inside $\MMM_{\cKur}$.
    \item[(1')] An arithmetic Breuil-Kisin-Fargues $G_K$-module $\MMM$ is in the essential image of $\MF^{wa}_{K}(\varphi,N)$ if and only if the initial group $I_K$ acts trivially on $\MMM_{\cKur}[\frac{1}{p}]$.
    \item[(2)] A Breuil-Kisin-Fargues $G_K$-module $\MMM$ is in the essential image of $\MF^{wa}_{K,\varphi}$ if and only if $(\MMM\otimes \overline{B})^{G_K}$ as a $K_0$-vector space has dimension equal to the rank of $\MMM$.
    \item[(2')] A Breuil-Kisin-Fargues $G_K$-module $\MMM$ is in the essential image of $\MF^{wa}_{K,\varphi}$ if and only if $(\MMM\otimes \overline{B})^{I_K}$ as a $\cKur$-vector space has dimension equal to the rank of $\MMM$.
\end{enumerate}
\end{thm}

\begin{rem}
In Theorem~\ref{essentialimages}, we have $(1)$ and $(1^\prime)$ (resp. $(2)$ and $(2^\prime)$) are equivalent due to unramified descent(cf. proof in Lemma~\ref{lemmaunrdescent}). Or this is from the well-known fact that a $p$-adic representation of $G_K$ is crystalline (resp. log-crystalline) if and only if when restricted on $I_K$ it is crystalline (resp. log-crystalline).
\end{rem}

For (2) in Theorem~\ref{essentialimages}, we need the following lemma.
\begin{lem}
Taking the $G_K$-invariant of the exact sequence
$$
0 \to \p \to \Ainf[\frac{1}{p}] \to \overline{B} \to 0,
$$
We will have $K_0\xhookrightarrow{} \overline{B}^{G_K}$. Let $\overline{B}^{G_K}_\mathrm{fin}$ be the sub-algebra consists of $v \in \overline{B}^{G_K}$ such that $\{\varphi^{n}(v)\}_{n\geq0}$ generates a finite dimensional $K_0$ vector space. We have $\overline{B}^{G_K}_\mathrm{fin}=K_0$.
\end{lem}

\begin{proof}
We have $\Ainf[\frac{1}{p}]^{G_K} = K_0$, i.e., $\Ainf[\frac{1}{p}]^{G_K}=W(\overline{k})[\frac{1}{p}]^{G_K}$. By the Teichm\"uller expansion of elements in $\p$, we will have 
$$\p^{G_K}=\p\cap W(\overline{k})[\frac{1}{p}]^{G_K}=0.$$ So we have $K_0\xhookrightarrow{} \overline{B}^{G_K}$.

Now $\overline{B}^{G_K}_\mathrm{fin}$ decomposes into isocrystals over $k$, in particular, it decomposes into subspace $(\overline{B}^{\varphi^h=p^d})^{G_K}$ for pairs $(h,d)$ with $d\in \Z$ and $h\in \N_{>0}$. Then by Lemma~\ref{Bdh}, we have
$$
(\overline{B}^{\varphi^h=p^d}) = (\Bcris^+)^{\varphi^h=p^d},
$$
which is also $G_K$-equivariant so the result follows from $(\Bcris^+)^{G_K}=K_0$.
\end{proof}

\begin{rem}
We don't know whether $\overline{B}^{G_K}=K_0$ at this time, but as in the statement of (2) in Theorem~\ref{essentialimages}, $(\MMM\otimes \overline{B})^{G_K}$ can be always viewed as a module over $\overline{B}^{G_K}_\mathrm{fin}=K_0$.
\end{rem}

\begin{proof}(of (1') Theorem~\ref{essentialimages}) 
If (the isogeny class of) an arithemtic BKF $G_K$-modules $\MMM$ corresponds to $D \in \MF^{wa}_{K}(\varphi,N,G_K)$ that also corresponds to a $G_K$-equivariant modification $\E_0\xdashrightarrow[]{} \E_1$, then the recall that as in \S 2.2, we show the $G_K$-action on $\E_1$ is given by $\alpha\wedge \E(D,\varphi,G_K))$. And recall $\alpha$ corresponds to the cocycle 
\begin{equation}\label{actionofG}
(\exp(-\log_{\underline{\varpi},g}))_g.
\end{equation}
On the other hand, we know by Proposition~\ref{BcrisofBKF}, $\E_1$ corresponds to $\MMM_\cris$. So there is a $G_K$-equivariant isomorphism
$$
D = (D\otimes \Bcris^+)\otimes_{\Bcris^+} \cKur = \MMM\otimes_{\Ainf} \cKur
$$
and $D \in \MF^{wa}_{K}(\varphi,N)$ if and only if $I_K$ acts trivially on $D$. So by equation~\ref{actionofG}, it is enough to show 
$$
\exp(-\log_{\underline{\varpi},g}N))\equiv 1 \mod{W(\m^\flat)\Bcris^+}.
$$
Recall we have $\log_{\underline{\varpi},g}$ is a multiple of $t$, so the result follows from $t \in W(\m^\flat)\Bcris^+$, which is directly from the definition of $t$.
\end{proof}

\begin{rem}\label{remark00}
If we assume an arithmetic BKF $G_K$-module $\MMM$ satisfies the condition (2) in Theorem~\ref{essentialimages}, then a similar argument will show that $\MMM$ corresponds to an object in $\MF^{wa}_{K,\varphi}$. However, we will prove something more substantial, i.e., condition (2) in Theorem~\ref{essentialimages} implies $\MMM$ being arithmetic!
\end{rem}

\begin{lem}\label{GKsection}
For $\MMM\in \BKF(G_K)$, let $M_{\cris}=\MMM\otimes \Bcris^+$, and $M_{\overline{B}}=\MMM\otimes_{\Ainf} \overline{B}$. Let $L_0=W(l)[\frac{1}{p}]$ with $k\subset l \subset \overline{k}$. For any $\varphi$-stable $L_0$-vector space $V\subset M_{\overline{B}}$ such that $V \otimes_{L_0} \overline{B} = M_{\overline{B}}$. Assume $V$ is also $G_K$-stable, then the $\varphi$-equivariant section $s: V\to M_\cris$ defined in Lemma \ref{lemmasection} is $G_K$-equivariant. 
\end{lem}

\begin{proof}
This is a slight generalization of \cite[Lemma 3.15]{Ozeki}. We use the same reduction as in Lemma \ref{lemmasection} and keep the notions. Fix $g\in G_K$ and let $D\in M_d(\Ainf)$ and $D_0\in M_d(\O_{L_0})$ be the matrix defined by 
$$
g(\overline{e}_1,\ldots,\overline{e}_d)=(\overline{e}_1,\ldots,\overline{e}_d)D_0
$$
and 
$$
g({e}_1,\ldots,{e}_d)=({e}_1,\ldots,{e}_d)D.
$$
We have $D-D_0=[\varpi']W$ with pseudo uniformizer $\varpi'\in C$ and $W\in M_d(\Ainf)$. It is enough to check 
$$
Dg(Y)=YD_0
$$
Since the $G_K$-action commutes with $\varphi$, we will have
$$
Dg(A)=A\varphi(D)\quad \text{ and }  \quad D_0g(A_0)=A_0\varphi(D_0).
$$

For $n\in \N_{>0}$, we have 
\begin{IEEEeqnarray*}{+rCl+x*}
Dg(Y_n) & = & Dg(A\varphi(A)\cdots \varphi^{n-1}(A)\varphi^{n-1}(A_0^{-1})\cdots \varphi(A_0^{-1})A_0^{-1}) \\
        & = & Dg(A)g(\varphi(A))\cdots g(\varphi^{n-1}(A))g(\varphi^{n-1}(A_0^{-1}))\cdots g(\varphi(A_0^{-1}))g(A_0^{-1}) \\
        & = & A\varphi(A)\cdots \varphi^{n-1}(A)\varphi^{n}(D)g(\varphi^{n-1}(A_0^{-1}))\cdots g(\varphi(A_0^{-1}))g(A_0^{-1}) \\
        & = & A\varphi(A)\cdots \varphi^{n-1}(A)\varphi^{n}(D_0+[\varpi']W)g(\varphi^{n-1}(A_0^{-1}))\cdots g(\varphi(A_0^{-1}))g(A_0^{-1}) \\
        & = & A\varphi(A)\cdots \varphi^{n-1}(A)g(\varphi^{n-1}(A_0^{-1}))\cdots g(\varphi(A_0^{-1}))g(A_0^{-1})D_0 \\
        &   & + \frac{[\varpi']^{p^n}}{p^{na}}A\varphi(A)\cdots \varphi^{n-1}(A)\varphi^{n}(W)g(\varphi^{n-1}(B_0))\cdots g(\varphi(B_0))g(B_0).
\end{IEEEeqnarray*}
And the result follows from that $\frac{[\varpi']^{p^n}}{p^{na}}$ converges to 0 in $B^+$.
\end{proof}

\begin{proof}(of (2) in Theorem~\ref{essentialimages})
By Lemma~\ref{GKsection}, we have the $G_K$-invariant section inside $(\MMM\otimes \overline{B})^{G_K}$ lifts to a $G_K$-invariant section inside $\MMM\otimes \Bcris^+$, so $(3)$ implies there is a $G_K$-fixed basis inside $\MMM\otimes \Bcris^+$. Base change along $\theta:\Bcris^+\to \Cp$, we have $\MMM\otimes \Cp$ has a $G_K$-invariant basis, i.e., $\MMM$ is arithmetic. 

Now as in Remark~\ref{remark00}, the rest is similar to the proof of (1). We give another way of proving by showing the $B_e$-part of $\E_1$ in the $G_K$-equivariant modification $\E_0\xdashrightarrow[]{} \E_1$ corresponds to $\MMM$ is crystalline in the sense that it is equal to $\mathscr{V}_{\log}(D)$ for some $\varphi$-module $D$ under the functor $\mathscr{V}_{\log}$ we defined in Proposition~\ref{Bereps}, see also Definition 10.2.13 in \cite{FF} for the notion of $B_e$-representations being crystalline. We have the $B_e$-representation corresponds to $\E_1$ is 
$$
M_e=(\MMM\otimes \Bcris)^{\varphi=1}.
$$
And by Proposition 10.2.12, $loc. cit.$, it is crystalline if and only if 
$$
\dim_{K_0}(M_e\otimes_{B_e} \Bcris)^{G_K}=\mathrm{rank}_{\Ainf} \MMM.
$$
But we have $D=(\MMM\otimes \overline{B})^{G_K}\subset(M_e\otimes_{B_e} \Bcris)^{G_K}$, so $M_e$ is crystalline. 

On the other hand, from the construction of the Fargues-Fontaine functor, when $N=0$, we have a $G_K$-equivariant isomorphism $\MMM_\cris\simeq D\otimes \Bcris^+$ and $D\otimes \Bcris^+$ is equipped with the diagonal action. When $\MMM$ corresponds to a crystalline representation, $D=D_\cris$ has a $K_0$-basis fixed by $G_K$. 
\end{proof}

\begin{rem}\label{remarkonBe}
 A key idea of Berger's $B_e$-representation theory is that for an $\O_X$-representation $\E$ of $G_K$, its $B_e$ part will determine the $(\varphi,N)$-module structure. This translates into to our theory says that for arithmetic Breuil-Kisin-Fargues modules $\MMM$, $G_K$-action on $\MMM\otimes \Bcris^+$ determines the $(\varphi,N)$-module structure, in particular, it tells $p$-adic Hodge property of $T(\MMM)$ being log-crystalline or crystalline.
\end{rem}

\begin{rem}
Since $T(\MMM)$ being log-crystalline or crystalline is just determined by $V(\MMM)$, it is easy to deduce the following ``integral version" of Theorem~\ref{essentialimages}.
\end{rem}

\begin{thm}\label{intversionofmaintheorem}
Let $\Rep_{\Z_p}^{dR}(G_K)$ (resp. $\Rep_{\Z_p}^{lcr}(G_K)$, resp. $\Rep_{\Z_p}^{cris}(G_K)$) be the category of de Rham (resp. log-crystalline, resp. crystalline) representations of $G_K$ over a $\Z_p$ lattice, then
\begin{enumerate}
    \item[(1)] There is a equivalence of $\Rep_{\Z_p}^{dR}(G_K)$ with the category of arithmetic BKF $G_K$-modules.
    \item[(2)] The essential image of $\Rep_{\Z_p}^{lcr}(G_K)$ of the functor in (1) are the arithmetic BKF modules such that there is a $G_K$-fixed basis in $\MMM_{\cKur}$.
    \item[(2')] The essential image of $\Rep_{\Z_p}^{lcr}(G_K)$ of the functor in (1) are the arithmetic BKF modules such that the $I_K$-action on $\MMM_{\cKur}$ is trivial.
    \item[(3)] The essential image of $\Rep_{\Z_p}^{cris}(G_K)$ of the functor in (1) are BKF $G_K$-modules such that $(\MMM\otimes \overline{B})^{G_K}$ as a $K_0$-vector space has dimension equal to the rank of $\MMM$.
    \item[(3')] The essential image of $\Rep_{\Z_p}^{cris}(G_K)$ of the functor in (1) are BKF $G_K$-modules such that $(\MMM\otimes \overline{B})^{I_K}$ as a $\cKur$-vector space has dimension equal to the rank of $\MMM$.
\end{enumerate}
\end{thm}

\section{Comparison of integral $p$-adic Hodge theories}\label{sectioncompat}
In this section we will compare our theory with some other existing theories in integral $p$-adic Hodge theory and prove certain compatibility results. First, we have the following rigidity result of arithmetic Breuil-Kisin-Fargues modules.
\begin{lem}(rigidity of arithmetic Breuil-Kisin-Fargues $G_K$-modules)\label{rigidityofBKF}
For any two arithmetic Breuil-Kisin-Fargues $G_K$-modules $\MMM_1$ and $\MMM_2$, if $T(\MMM_1)\simeq T(\MMM_2)$, then $\MMM_1\simeq\MMM_2$.
\end{lem}
\begin{proof}
By Fargues classification theorem (cf. Theorem~\ref{Farguesclassification}), it is enough to show there is a $G_K$-equivariant isomorphism
$$
\MMM_1\otimes\BdR^+ \simeq \MMM_2\otimes\BdR^+.
$$
But we have known by definition, both $\MMM_1\otimes\BdR^+$ and $\MMM_2\otimes\BdR^+$ are $\BdR^+$-flat representations inside $T(\MMM_1)\otimes\BdR$, so both of them are isomorphic to $D_{\dR}(T(\MMM))\otimes_K \BdR^+$.
\end{proof}

\subsection{Wach modules and crystalline representations}
 Throughout this section, we assume that $K = K_0$ is unramified for simplicity. We will discuss general $K$ in Remark~\ref{WachKisinRen}. We will construct a functor from Wach modules over $K$ to $\BKF^a(G_K)$. First recall the definition of Wach modules after Berger. Let $\{\epsilon_i\}_{i\geq 0}$ be a compatible system of $p^n$-th roots of 1, i.e., $\epsilon_0=1$ and $\epsilon_{i+1}^p=\epsilon_{i}$ for all $i$. $\{\epsilon_i\}_{i\geq 0}$ defines an element $\underline{\epsilon}\in \O_C$, let $\mu=[\underline{\epsilon}]-1\in \Ainf$. Let $\SSS_\mu=W(k)[[\mu]]$ as a subring of $\Ainf$ which is stable under $\varphi$ on $\Ainf$. Let $q=\frac{\varphi(\mu)}{\mu}\in \SSS_\mu$. By \cite[\S 3.3]{BMS1}, we have $q$ generates the kernel of $\theta\circ\varphi^{-1}$, i.e., we have $(q)=(\tilde{\xi})$ in $\Ainf$. Let $K_{p^\infty}=\cup_i K(\epsilon_i)$, and let $\Gamma=\Gal(K_{p^\infty}/K)$, note that $\Gamma$ acts on $\SSS_\mu$.
\begin{defn}
A (finite free) Wach module $\MM$ over $K$ is a finite free $\SSS_{\mu}$ module together with
$$
\varphi_{\MM_\mu}: \SSS_\mu\otimes_{\varphi,\SSS_\mu} \MM_\mu \to \MM_\mu
$$
such that the cokernal is killed by a power of $q$. $\MM_\mu$ also equipped with a semilinear $\Gamma$-action on $\MM_\mu$ commutes with $\varphi_{\MM_\mu}$, satisfies $\Gamma$ acts trivially on $\overline{\MM_\mu}:=\MM_\mu\otimes_{\SSS_\mu}W(k)$. 
\end{defn}

And we have the following theorem of Berger.
\begin{thm}(\cite{Bergeronwachmod})
There is an equivalence of $G_K$-stable $\Z_p$-lattices in crystalline representations with non-negative Hodge–Tate weights and the category of finite free Wach modules. And satisfies if $\MM_{\mu}$ corresponds to a crystalline representation $T$, then 
$$
(\MM_{\mu}\otimes_{\SSS_{\mu}}W(C))^{\varphi=1}\simeq T.
$$
\end{thm}

\begin{rem}
Here we use the covariant version of the theorem, note that we still get non-negative Hodge–Tate weights, since in our convention~\ref{conventionHTw} we let the cyclotomic character to have Hodge-Tate weight $-1$.
\end{rem}

If $\MM$ is a Wach module, let define $\MMM_{\mu}=\MM_\mu\otimes_{\SSS_\mu} \Ainf$, $\MMM_{\mu}$ has a $\varphi$-structure defined by $\varphi_{\MM_\mu}\otimes\varphi_{\Ainf}$ and a semilinear $G_K$-action.

\begin{lem}\label{WachtoBKF}
$\MMM_{\mu}$ with the $G_K$-action defined as above is an arithmetic Breuil-Kisin-Fargues $G_K$-module, moreover, it is actually crystalline in the sense that it satisfies the condition (3) in Theorem~\ref{intversionofmaintheorem}.
\end{lem}
\begin{proof}
Firstly, we have $(q)=(\tilde{\xi})$, so $\MMM_{\mu}$ together with $G_K$-action is a Breuil-Kisin-Fargues $G_K$-module. We claim $(\MMM_{\mu}\otimes\Bcris^+)^{G_K}$ has a basis over $K_0$ of full rank. This is from the fact there is a $G_K$-invariant section of $\overline{\MM}$ inside $\MMM_{\mu}\otimes\Bcris^+$ by \cite[Lemma 2.2.2]{KisinRen}. 
\end{proof}

\begin{rem}\label{WachKisinRen}
In the proof of Lemma~\ref{WachtoBKF}, when $K=K_0$, we should also have $\MMM_{\mu}$ is arithmetic for Kisin-Ren's generalization of Wach modules (c.f. \cite{KisinRen}), where the cyclotomic tower is replaced by Lubin-Tate tower. For general ramified $K$, there should be a paralleled theory for arithmetic Breuil-Kisin-Fargues $G_K$-modules over $A_{\mathrm{inf},K} := \Ainf\otimes_{W(k)}\O_K$ that relates to modifications of vector bundles over the Fargues-Fontaine curve defined using perfectoid field $C$ and the field $K$ (rather than $\Q_p$). We will discuss this in further work.
\end{rem}

\subsection{Kisin modules, (\texorpdfstring{$\varphi,\hat{G}$}{f,g})-modules}

For general ramified $K$. We will use $\vec{\varpi}:=\{\varpi_n\}$ to denote a compatible system of $p^n$-th roots of a uniformizer $\varpi$ of $\O_K$, i.e., $\varpi_0=\varpi$ and $\varpi_{i}=\varpi_{i+1}^p$ for all $i$. We define $K_{\infty}=\cup_{n=1}^{\infty} K(\varpi_n)$, we will also write it as $K_{\infty,\vec{\varpi}}$ when we want to emphasize the choice of $\vec{\varpi}$. The compatible system $\{\varpi_n\}$ also defines an element $\underline{\varpi}$ in $\O_C$. View $\SSS=W(k)[[u]]$ as a sub-$W(k)$-algebra of $\Ainf$ determined by $u \mapsto [\underline{\varpi}]$. We will also use $\SSS_{\vec{\varpi}}$ to emphasize the choice of $\vec{\varpi}$. One can check $\varphi_{\Ainf}(u)=u^p$, so in particular $\SSS$ in stable under $\varphi_{\Ainf}$, let $\varphi_\SSS=\varphi_{\Ainf}|_\SSS$. We also have $G_{K_\infty}$ fix $u$ so $G_{K_\infty}$ acts trivially on $\SSS$. Let $E(u)\in W(k)[u]$ be a minimal polynomial of $\varpi$ over $K_0$. Also let $K_{\infty,p^\infty}=\cup_iK(\varpi_i,\epsilon_i)$, i.e., $L$ is the normalization of $K_{\infty}$, define $\widehat{G}=\Gal(K_{\infty,p^\infty}/K)$.

\begin{defn}
A (finite free) Kisin module is finite free $\SSS$ module $\MM$ together with
$$
\varphi_{\MM}: \SSS\otimes_{\varphi,\SSS} \MM \to \MM
$$
such that the cokernel of $\varphi_{\MM}$ is killed by a power of $E(u)$.
\end{defn}

One can check $E(u)$ generates kernel of $\theta$, i.e., $(\xi)=(E(u))$ in $\Ainf$, so we have
\begin{lem}
For a Kisin module $\MM$, $\MM\otimes_{\SSS,\varphi}\Ainf$ is a Breuil-Kisin-Fargues module. 
\end{lem}

Let $T$ be a log-crystalline representation of $G_K$ over $\Z_p$ with non-negative Hodge-Tate weights, then Kisin in \cite{KisinFcrystal} can show $T$ is of finite $E(u)$-height in the sense that there is Kisin module 
$$
\MM \subset (T\otimes_{\Z_p} W(C))^{G_{K_\infty}}
$$
that is $\varphi$-stable and spans $T\otimes_{\Z_p} W(C)$ as a $W(C)$-module. This Kisin module is also uniquely determined by $T = (\MM\otimes_{\SSS}W(C))^{\varphi=1}$ as a representation of $G_{K_{\infty}}$. So we can write $\MM=\MM(T)$ and let $\MMM(T)=\MM(T)\otimes_{\SSS,\varphi}\Ainf$ be the corresponded Breuil-Kisin-Fargues module which carries a natural $G_{K_{\infty}}$-semilinear action that commutes with $\varphi_{\MMM}$. We claim that there is a unique way to extend this to a $G_K$-semilinear action so that $\MMM(T)$ is an arithmetic Breuil-Kisin-Fargues $G_K$-module. Actually, we will prove a more general result for $T$ being potentially log-crystalline representation of $G_K$ over $\Z_p$ with non-negative Hodge-Tate weights, we choose any finite Galois extension $L/K$ such that $T|_{G_L}$ is log-crystalline. Then as we have discussed, $T_L:=T|_{G_L}$ is of finite $E$-height for a choice of $p^n$-th roots of a uniformizer $\varpi_L$ of $\O_L$ and let $\MM(T_L)$ (resp. $\MMM(T_L)$) be the corresponded Kisin module (resp. Breuil-Kisin-Fargues module with $G_{L_\infty}$-action).

\begin{thm}
Up to isomorphic, there is a unique way to extends the $G_{L_{\infty}}$-semilinear action on $\MMM(T):=\MMM(T_L)$ to an action of $G_K$ so that $\MMM(T)$ is an arithmetic Breuil-Kisin-Fargues $G_K$-module and satisfies a $G_K$-equivariant isomorphism
$$
(\MMM(T)\otimes W(C))^{\varphi=1}=T.
$$
\end{thm}

\begin{proof}
The proof is an explicit comparing the construction of Kisin module and our construction of arithmetic Breuil-Kisin-Fargues modules. We need to compute $(T,\Xi)$ and its $G_{L_\infty}$-action for $\MMM(T_L)$. 

We give a brief review of the construction of Kisin module from log-crystalline representation: let $T$ is a potentially log-crystalline representation of $G_K$ over $\Z_p$, $L/K$ be a finite Galois extension such that $T|_{G_L}$ becomes log-crystalline, and let $L_0=W(k_L)[\frac{1}{p}]$ and define $D'=(T\otimes \Bst)^{G_L}$ as the filtered $(\varphi,N,G_K)$-module associated with $T\otimes \Q_p$. Then we obtain a filtered $(\varphi,N)$-module $D'$ over $L$ or $(D',\varphi,N,\Fi^\bullet)$ by forgetting the $G_K$-action. $D'$ corresponds to the log-crystalline representation $T\otimes \Q_p|_{G_L}$. Now let $\O$ be the ring of rigid analytic functions over the open unit disc over $L_0$ in the variable $u$. Let $\SSS=W(k_L)[[u]]$, then one has $\SSS[\frac{1}{p}]\subset \O$ and there is a $\varphi_\O$ extending $\varphi_\SSS$. Fix $(\varpi_{L,n})$ any choice of compatible system of $p^n$-th roots of a uniformizer $\varpi_{L,0}$ of $L$, then one can easily show that the the inclusion $\SSS[\frac{1}{p}]\to \Ainf[\frac{1}{p}]$ with $u \mapsto [(\overline{\varpi_{L,n}})]$ extends to an inclusion $\O\to B^+$. Geometrically, $\O$ (resp. $B^+$) is the locus $\{p\neq 0\}$ of $\Spa(\SSS)$ (resp. $\Spa(\Ainf)$), and restrict the covering map $\Spa(\Ainf)\to \Spa(\SSS)$ to these loci will give $\O\to B^+$. 

Roughly speaking, Kisin defined the $\SSS$ module by descending a $\varphi$-module $\M(D')$ over $\O$ using the theory of slope of Kedlaya. In particular, we have $\MM\otimes\O=\M(D')$. And a theorem of Fontaine says that the ways of descent $\M(D')$ to $\MM$ are canonically corresponded with $G_{L_{\infty}}$-stable $\Z_p$-lattices in $T\otimes \Q_p$, where $L_{\infty}=\cup_{n=1}^{\infty} L(\varpi_{L,n})$. Then Kisin define $\MM$ to be the $\SSS$-module descents $\M(D')$ using the lattice $T|_{G_{L_{\infty}}}$. So we have
$$
T(\MMM(T_L))=(\MM\otimes_\SSS W(C))^{\varphi=1}=T|_{G_{L_{\infty}}}.
$$
For $\Xi=\MMM(T_L)\otimes \BdR^+$, we need to review the construction of $\M(D')$ we mentioned above. For all $n\in \Z$ consider the composition:
\[
\begin{tikzcd}
\theta_n : \O \arrow[r] & B^+ \arrow[r,"\varphi^{-n}"] & B^+ \arrow[r, "\theta"] & \Cp
\end{tikzcd}
\]
and let $x_n$ be the closed points on the rigid open unit disc defined by $\theta_n$. And define $\O_{st}=\O[l_u]$, the $\O$-algebra generated by $l_u=\underline{\varpi}_L]$ inside $\BdR^+$. And extend the $\varphi$-action to $l_u$ by $\varphi(l_u)=p l_u$ and define a $\O$-derivation $N$ on $\O^+_{st}$ by letting $N(l_u)=1$. Given $(D',\varphi,N,\Fi^\bullet)$, Kisin defines $\M(D')$ as certain modification of the vector bundle $$(\O[l_u]\otimes_{L_0} D)^{N=0}$$ over $\O$ along stalks at $x_n$ for $n\geq 0$. However to compute $\Xi$, we only need to the stalk at $x_{-1}$. and there is a natural isomorphism\cite[Proposition 1.2.8.]{KisinFcrystal}:
\[
\begin{tikzcd}
(\O[l_u]\otimes_{L_0} D')^{N=0} = (L_0[l_u]\otimes_{L_0} D')^{N=0}\otimes_{L_0}\O \arrow[r, "\eta\otimes \mathrm{id}"] & D'\otimes_{L_0}\O.
\end{tikzcd}
\]
This tells us $\Xi=(\M(D')\otimes_{\O}B^+)\otimes_{B^+,\varphi} \BdR^+$ is isomorphic to
$$
(D'\otimes_{L_0,\varphi} L_0)\otimes_{L_0}\BdR^+ \cong D'\otimes_{L_0} \BdR^+
$$ 
with $G_{L_\infty}$ acts trivially on $D'$ by construction. 

To finish the proof, we just need to show there is a unique way to extend the $G_{L_\infty}$-action on $(T_{L_\infty},D'\otimes \BdR^+)$ to $G_K$ such that on $T$ it is the original potentially log-crystalline representation, and on $\Xi=D'\otimes \BdR^+$ is $\BdR^+$-flat. But we have if we want $\Xi\in T\otimes \BdR$ is $G_K$ stable and flat, by Proposition~\ref{BdR+flat}, it has to equal to 
$$
D_{\dR}(T)\otimes_{K} \BdR^+=(D'\otimes L)^{G_K}\otimes \BdR^+.
$$
\end{proof}

We also recall the following definition of $(\varphi,\widehat{G})$-modules of Liu.
\begin{defn}
Let $\mathrm{Mod}^{\varphi,\widehat{G}}_{\SSS,\widehat{R}}$ be the category of triples $(\MM,\varphi_{\MM},\widehat{G})$ called $(\varphi,\widehat{G})$-modules, where
\begin{enumerate}
    \item $(\MM,\varphi_\MM)$ is a Kisin module;
    \item $\widehat{G}$ is a continuous $\widehat{R}$-semilinear $\widehat{R}$-action on $\widehat{\MM}:=\MM\otimes_{\SSS,\varphi}\widehat{R}$;
    \item $\widehat{G}$ commutes with $\varphi_{\widehat{\MM}}$;
    \item regarding $\MM$ as a $\varphi(\SSS)$-submodule of $\widehat{\MM}$, we have $\MM\subset \widehat{\MM}^{\Gal(K_{\infty,p^\infty}/K_{\infty})}$;
    \item $\widehat{G}$ acts trivially on $\widehat{\MM}\otimes_{\widehat{R}}W(k)$.
\end{enumerate}
\end{defn}

Here $\widehat{R}$ is a subring of $\Ainf$. We will not give the explicit definition of $\widehat{R}$, we just list two properties we need in our applications: 
\begin{enumerate}
    \item $\widehat{R}\subset \Ainf^{G_{K_{\infty,p^\infty}}}$;
    \item the image of $\widehat{R} \xhookrightarrow{} \Ainf \xrightarrow{\theta}$ is $K$.
\end{enumerate}

The main result in \cite{liu-notelattice} is
\begin{thm}
There is an equivalence of log-crystalline representations of non-negative Hodge-Tate weights with the category of $(\varphi,\hat{G})$-modules. And satisfies if $(\MM,\varphi_{\MM},\widehat{G})$ corresponds to a log-crystalline representation $T$, then 
$$
(\widehat{\MM}\otimes_{\SSS}W(C))^{\varphi=1}\simeq T.
$$
\end{thm}

\begin{lem}\label{ghatmodule}
$(\MM,\varphi_{\MM},\widehat{G})$ be a $(\varphi,\widehat{G})$-module, let $\MMM=\widehat{\MM}\otimes_{\widehat{R}}\Ainf$, then $\MMM$ is a Breuil-Kisin-Fargues $G_K$-module where the $G_K$-action comes from the $\widehat{G}$-action on $\widehat{\MM}$, moreover, $\MMM$ is arithmetic.
\end{lem}
\begin{proof}
We have $\MMM=\widehat{\MM}\otimes_{\widehat{R}}\Ainf=\MM\otimes_{\SSS,\varphi}\Ainf$ is a Breuil-Kisin-Fargues module. And we have $\widehat{\MM}\otimes_{\widehat{R},\theta}K$ is a $K$ vector space and
$$
\MMM\otimes \Cp=(\widehat{\MM}\otimes_{\widehat{R},\theta}K)\otimes_K \Cp.
$$
Moreover, use the definition, we have $\widehat{\MM}\otimes_{\widehat{R},\theta}K$ has a $G_K$ linear action and there is a basis fixed by $G_{K_\infty}$ coming from a basis of $\MM$. The result is from the following lemma on Kummer theory and Galois descent.
\end{proof}

\begin{lem}
All closed normal subgroup of $G_K$ containing $G_{K_\infty}$ are open.
\end{lem}
\begin{proof}
This is from the fact $K$ can only contain finitely mainly $p^n$-th roots of $1$.
\end{proof}

So now as we mentioned in the beginning of \S~\ref{sectioncompat}, we have the following compatibility results for Kisin modules defined using different choice of uniformizer and Kummer tower, and the compatibility of Kisin modules, Wach modules and Kisin-Ren's modules when $K=K_0$ and $T$ is crystalline.
\begin{thm}\label{compatible}
Let $T$ be a log-crystalline representations with non-negative Hodge-Tate weights, for different choice of $\vec{\varpi}$, the arithmetic Breuil-Kisin-Fargues $G_K$-modules $\MMM$ from Lemma~\ref{ghatmodule} are all isomorphic to the arithmetic Breuil-Kisin-Fargues modules $\MMM(T)$ defined in Theorem~\ref{intversionofmaintheorem}, in particular, they are all isomorphic to each other. 

Moreover, if $T$ is crystalline then $\MMM(T)$ is isomorphic to the arithmetic Breuil-Kisin-Fargues $G_K$-modules defined from Wach modules as in Lemma~\ref{WachtoBKF}.
\end{thm}

\begin{rem}
We want to note that the above result also been proved in \cite{liu2013compatibility}.
\end{rem}

\subsection{Prismatic $F$-crystals are arithmetic}\label{subsec:prismaticarearithmetic}
In this subsection, we will review the theory of prismatic $F$-crystals over $\O_K$ and its logarithmic variation without carefully going through every definition, all the notations we are using is from \cite{BS19}, \cite{Bhatt-Scholze_prismaticFcrystal} and \cite{DuLiu_prismaticphihatG}. 

First recall over $\O_K$, we can define its absolute prismatic site $(\O_K)_\Prism$ (resp. logarithmic prismatic site $(\O_K)_{\Prism_{\log}}$) with structure sheaf $\O_\Prism$ that admits an endomorphism $\varphi_\Prism$ together with the Hodge-Tate ideal sheaf $\mathcal{I}_\Prism \subset \O_\Prism$. 

\begin{defn}
Let $\ast\in\{\emptyset,\log\}$, a prismatic $F$-crystal over $(\O_K)_{\Prism_{\ast}}$ is a crystal $\MM_\Prism$ in finite locally free $\O_\Prism$-modules with an isomorphism
$$
(\varphi^\ast\MM_\Prism)[\frac{1}{\mathcal{I_\Prism}}] \simeq \MM_\Prism[\frac{1}{\mathcal{I_\Prism}}].
$$
\end{defn}

An important result is the following
\begin{thm}(\cite[Theorem 1.2]{Bhatt-Scholze_prismaticFcrystal}, \cite[Theorem 1.0.2]{DuLiu_prismaticphihatG})
There is an equivalence between the category of prismatic $F$-crystal over $(\O_K)_{\Prism}$ (resp. prismatic $F$-crystals over $(\O_K)_{\Prism_{\log}}$) with the category lattices in crystalline (resp. log-crystalline) representations of $G_K$.
\end{thm}

\begin{rem}\label{rem:remarkonetalerealization}
The functor in the above theorem can be explicitly written in the following way. Let $(\Ainf, \ker\theta)$ be the infinitesimal prism of Fontaine, then $(\Ainf, \ker\theta)$ admits an action of $G_K$ as an object in $(\O_K)_\Prism$. Let $\MM_\Prism$ be a prismatic $F$-crystal, then $\MM_\Prism[1/[\underline{\varpi}^\flat]]^\wedge_p$ is a finite free \'etale $\varphi$-module over $W(C)$-module, and we define $T(\MM_\Prism)$ to be finite free $\Z_p$-module $(\MM_\Prism[1/[\varpi^\flat]]^\wedge_p)^{\varphi=1}$ as representation of $G_K$.
\end{rem}

We have the following result.
\begin{thm}\label{prismaticarearithmetic}
Let $(\Ainf,\varphi(\ker\theta))$ be the twisted infinitesimal prism equipped with the $G_K$-action, then for any (logarithmic) prismatic $F$-crystal $\MM_\Prism$, we have $\MM_\Prism(\Ainf,\varphi(\ker\theta))$ is an arithmatic BKF-module.
\end{thm}
\begin{proof}
First from the definition of prismatic $F$-crystal, $\MM_\Prism(\Ainf,\varphi(\ker\theta))$ is a BKF-module. To check $\MM_\Prism(\Ainf,\varphi(\ker\theta))$ is arithmetic, one has by \cite[\S~3.3]{DuLiu_prismaticphihatG}, $\MM_\Prism(\Ainf,\varphi(\ker\theta))$ comes from the base change from a $(\varphi,\hat{G})$-module. So by Theorem~\ref{compatible}, we have it is arithmetic.
\end{proof}

Theorem~\ref{prismaticarearithmetic} together with Remark~\ref{rem:remarkonetalerealization} enables us to compare $\MM_\Prism(\Ainf,\varphi(\ker\theta))$ with other integral $p$-adic Hodge theories at $\Ainf$-level.

\subsection{Breuil-Kisin-Fargues \texorpdfstring{$G_K$}{GK}-modules admit all descents}\label{subsectionalldescent}

\begin{defn}\cite[F.7. Definition]{emertongee2020moduli}\label{alldescent}
Let $\MMM$ be a Breuil-Kisin-Fargues $G_K$-module. Then we say that $\MMM$ admits all descents over $K$ if the following conditions hold.
\begin{enumerate}
\item For any choice $\varpi$ of uniformaizer of $\O_K$ and any compatible system $\vec{\varpi}=(\varpi_n)$ of $p^n$-th roots of $\varpi$, there is a Breuil-Kisin module $\MM_{\vec{\varpi}}$ defined using $\vec{\varpi}$ such that $\MM_{\vec{\varpi}}\otimes_{\SSS,\varphi}\Ainf$ is isomorphic to $\MMM$ and $\MM_{\vec{\varpi}}$ is fixed by $G_{K_{\vec{\varpi},\infty}}$ under the above isomorphism, where $K_{\vec{\varpi},\infty}=\cup_n K(\varpi_n)$;
\item $\MM_{\vec{\varpi}}\otimes_{\SSS,\varphi}(\SSS/E(u_{\vec{\varpi}})\SSS)$ is independent of the choice of $\vec{\varpi}$ as a $\O_K$-submodule of $\MMM\otimes_{\Ainf}\O_{\Cp}$;
\item Let $u_{\vec{\varpi}}=[(\overline{\varpi_n})]$, then $\MM_{\vec{\varpi}}\otimes_{\SSS,\varphi}(\SSS/u_{\vec{\varpi}}\SSS)$ is independent of the choice of $\vec{\varpi}$ as a $W(k)$-submodule of $\MMM\otimes_{\Ainf}W(\overline{k})$.
\end{enumerate}
\end{defn}

\begin{rem}
Theorem~\ref{compatible} will implies if $T$ is a log-crystalline representation of $G_K$, and let $\MMM(T)$ be arithmetic Breuil-Kisin-Fargues $G_K$-module corresponds to $T$ under Theorem~\ref{intversionofmaintheorem}, then $\MMM(T)$ admits all descents over $K$.
\end{rem}

The following result is first proved by Gee-Liu in \cite{emertongee2020moduli}, which can be regarded as an inverse of Theorem~\ref{compatible}. 
\begin{prop}\label{p:7}
Let $\MMM$ be a Breuil-Kisin-Fargues $G_K$-module, and assume $\MMM$ admits all descents over $K$, then $\MMM$ is arithmetic and satisfies the condition $(2)$ in Theorem~\ref{intversionofmaintheorem}, i.e., the inertia subgroup $I_K$ of $G_K$ acts trivially on $\overline{\MMM}= \MMM\otimes_{\Ainf}W(\overline{k})$. In particular, $T(\MMM)$ is log-crystalline.
\end{prop}
\begin{proof}
The proof is bases on \cite[F.15]{emertongee2020moduli} on Kummer theory, which will implies the closed subgroup generated by $\{K_{\vec{\varpi},\infty}\}_{\vec{\varpi}}$ is $G_K$.

So $(2)$ in Definition~\ref{alldescent} will imply $\MMM\otimes\Cp$ has a $G_K$ fixed basis, i.e., $\MMM$ is arithmetic. And $(3)$ in Definition~\ref{alldescent} will imply that $\MMM\otimes\cKur$ also has a $G_K$ fixed basis, so the the inertia subgroup $I_K$ acts trivially on on this basis.
\end{proof}

\begin{rem}
In a recent work \cite{Gao2020breuilkisin}, he was able to show that if an arithmetic Breuil-Kisin-Fargues $G_K$-module $\MMM$ admits one descent, that is for one choice of compatible system $\vec{\varpi}=(\varpi_n)$ of $p^n$-th roots of a uniformizer $\varpi$ of $\O_K$, there is a Breuil-Kisin module $\MM_{\vec{\varpi}}$ defined using $\vec{\varpi}$ such that $\MM_{\vec{\varpi}}\otimes_{\SSS,\varphi}\Ainf$ is isomorphic to $\MMM$ and $\MM_{\vec{\varpi}}$ is fixed by $G_{K_{\vec{\varpi},\infty}}$ under the above isomorphism, then $\MMM$ is arithmetic. Actually, his work corrects a mistake in \cite{Caruso} and shows that if a $p$-adic representation is of finite $E$-height, then the representation is de Rham. We want to note this result is deep since for a representation with non-negative Hodge-Tate weights is of finite $E$-height that defined by a Kisin module
$$
\MM \xhookrightarrow{} T\otimes_{\Z_p} W(C)
$$
then we don't know a priori that $\MM\otimes_{\SSS,\varphi}\Ainf \xhookrightarrow{} T\otimes_{\Z_p} W(C)$ is stable $G_K$. 
\end{rem}

\subsection{Crystalline condition}
At last, we show how the ring $\overline{B}$ helps define a crystalline condition for arithmetic Breuil-Kisin-Fargues $G_K$-module $\MMM$ with descents to Kisin modules. We observe that for Kisin's $\SSS$
$$
\SSS \xrightarrow{\varphi} \Ainf \to \overline{B}
$$
uniquely factor through $\SSS \to W(k)$ since $u\in p$. So we slightly modify the $(3)$ in Definition~\ref{alldescent}
\begin{defn}\label{crystalline}
Let $\MMM$ be a Breuil-Kisin-Fargues $G_K$-module, and $\MMM$ admits all descents over $K$, we say it is crystalline if it also satisfies
\begin{enumerate}
\item[(3')] Let $u_{\vec{\varpi}}=[(\overline{\varpi_n})]$, then $\MM_{\vec{\varpi}}\otimes_{\SSS,\varphi}(\SSS/u_{\vec{\varpi}}\SSS)$ is independent of the choice of $\vec{\varpi}$ as a $W(k)$-submodule of $\MMM\otimes_{\Ainf}\overline{B}$.
\end{enumerate}
\end{defn}

A straightforward consequence from (3) in Theorem~\ref{intversionofmaintheorem} is
\begin{prop}\label{p:8}
Let $\MMM$ be a Breuil-Kisin-Fargues $G_K$-module and assume $\MMM$ satisfies $(3')$ in Definition~\ref{crystalline}, then $\MMM$ is arithmetic and satisfies the condition $(3)$ in Theorem~\ref{intversionofmaintheorem}. In particular, $T(\MMM)$ is crystalline.
\end{prop}

\begin{rem}
The crystalline conditions used in \cite{emertongee2020moduli} and \cite{Gao2020breuilkisin} will automatically fit into our Proposition~\ref{p:8}, for example, in \cite{emertongee2020moduli}, we have for $u=[\underline{\varpi}]$, $gu-u \in \p$ for all $g\in G_K$. 
\end{rem}

\bibliographystyle{amsplain}
\bibliography{library}

\def\cprime{$'$}
\providecommand{\bysame}{\leavevmode\hbox to3em{\hrulefill}\thinspace}
\providecommand{\MR}{\relax\ifhmode\unskip\space\fi MR }
\providecommand{\MRhref}[2]{%
  \href{http://www.ams.org/mathscinet-getitem?mr=#1}{#2}
}
\providecommand{\href}[2]{#2}
\begin{thebibliography}{10}

\bibitem{ansch_BKFCM}
Johannes Ansch\"utz, \emph{Breuil-{K}isin-{F}argues modules with complex
  multiplication}, Journal of the Institute of Mathematics of Jussieu (2020),
  1–50.

\bibitem{ALB}
Johannes Ansch{\"u}tz and Arthur-C{\'e}sar~Le Bras, \emph{{Prismatic
  Dieudonn{\'e} theory}}, 2020, arXiv:1907.10525.

\bibitem{Bergeronwachmod}
Laurent Berger, \emph{Limites de repr\'esentations cristallines}, Compos. Math.
  \textbf{140} (2004), no.~6, 1473--1498. \MR{2098398 (2006c:11138)}

\bibitem{BergerBpairs}
Laurent Berger, \emph{Construction de $(\varphi,{\Gamma})$-modules:
  repr\'esentations $p$-adiques et {$B$}-paires}, Algebra Number Theory
  \textbf{02} (2009), 91--120.

\bibitem{BMS1}
Bhargav Bhatt, Matthew Morrow, and Peter Scholze, \emph{Integral {$p$}-adic
  {H}odge theory}, Publ. Math. Inst. Hautes \'{E}tudes Sci. \textbf{128}
  (2018), 219--397. \MR{3905467}

\bibitem{BS19}
Bhargav Bhatt and Peter Scholze, \emph{{{Prisms and Prismatic Cohomology}}},
  2019, arXiv:1905.08229.

\bibitem{Bhatt-Scholze_prismaticFcrystal}
Bhargav Bhatt and Peter Scholze, \emph{Prismatic {$F$}-crystals and crystalline
  {G}alois representations}, 2021.

\bibitem{CMInotes}
Olivier Brinon and Brian Conrad, \emph{{CMI} summer school notes on $p$-adic
  hodge theory}, \url{https://math.stanford.edu/~conrad/papers/notes.pdf}.

\bibitem{Cais-Liu}
Bryden Cais and Tong Liu, \emph{On {$F$}-crystalline representations}, Doc.
  Math. \textbf{21} (2016), 223--270. \MR{3505135}

\bibitem{Caruso}
Xavier {Caruso}, \emph{{Repr{\'e}sentations galoisiennes {p}-adiques et
  ($\phi,\tau$)-modules}}, {Duke Mathematical Journal} \textbf{162} (2013),
  no.~13, 2525--2607.

\bibitem{prefaceofFF}
Pierre Colmez, \emph{La courbe de fargues et fontaine.}, Preface to \cite{FF}
  (2017).

\bibitem{CorIris}
Christophe Cornut and Macarena Peche~Irissarry, \emph{Harder--{N}arasimhan
  filtrations for breuil--kisin--fargues modules}, Annales Henri Lebesgue
  \textbf{2} (2019), 415--480 (en). \MR{4015914}

\bibitem{HDAinf}
Heng Du, \emph{$\mathbf{A}_{\text {inf}}$ has uncountable krull dimension},
  2020, arXiv:2002.10358.

\bibitem{DuLiu_prismaticphihatG}
Heng Du and Tong Liu, \emph{A prismatic approach to $(\varphi,
  \hat{G})$-modules and {$F$}-crystals}, 2021.

\bibitem{emertongee2020moduli}
Matthew Emerton and Toby Gee, \emph{{Moduli stacks of {\'e}tale
  ($\phi,\Gamma$)-modules and the existence of crystalline lifts}}, 2020,
  arXiv:1908.07185.

\bibitem{Faltingssimpson1}
Gerd Faltings, \emph{A {{\(p\)}}-adic {Simpson} correspondence}, Adv. Math.
  \textbf{198} (2005), no.~2, 847--862 (English).

\bibitem{Farguesconjecture}
Laurent Fargues, \emph{Quelques r\'{e}sultats et conjectures concernant la
  courbe}, Ast\'{e}risque (2015), no.~369, 325--374. \MR{3379639}

\bibitem{FF}
Laurent Fargues and J.-M Fontaine, \emph{Courbes et fibr\'es vectoriels en
  th\'eorie de hodge $p$-adique}, Ast\'erisque 406.

\bibitem{Gao2020breuilkisin}
Hui Gao, \emph{Breuil-{K}isin modules and integral $p$-adic hodge theory},
  Journal of the European Mathematical Society (2020).

\bibitem{How}
Sean Howe, \emph{Transcendence of the hodge-tate filtration}, Journal de
  Th\'eorie des Nombres de Bordeaux \textbf{30} (2018), no.~2, 671--680 (en).

\bibitem{KedlayaLiu1}
Kiran~S. Kedlaya and Ruochuan Liu, \emph{Relative {{\(p\)}}-adic {Hodge}
  theory: foundations}, Ast{\'e}risque, vol. 371, Paris: Soci{\'e}t{\'e}
  Math{\'e}matique de France (SMF), 2015 (English).

\bibitem{KisinFcrystal}
Mark Kisin, \emph{Crystalline representations and {$F$}-crystals}, Algebraic
  geometry and number theory, Progr. Math., vol. 253, Birkh\"auser Boston,
  Boston, MA, 2006, pp.~459--496. \MR{MR2263197 (2007j:11163)}

\bibitem{Kisin2adic}
\bysame, \emph{Modularity of 2-adic {B}arsotti-{T}ate representations}, Invent.
  Math. \textbf{178} (2009), no.~3, 587--634. \MR{2551765}

\bibitem{KisinRen}
Mark Kisin and Wei Ren, \emph{Galois representations and {L}ubin-{T}ate
  groups}, Doc. Math. \textbf{14} (2009), 441--461. \MR{2565906 (2011d:11122)}

\bibitem{LLAinf}
J~Lang and J~Ludwig, \emph{Ainf is infinite dimensional}, Journal of the
  Institute of Mathematics of Jussieu (2020).

\bibitem{liu-notelattice}
Tong Liu, \emph{A note on lattices in semi-stable representations}, Math. Ann.
  \textbf{346} (2010), no.~1, 117--138. \MR{2558890}

\bibitem{liu2013compatibility}
Tong Liu, \emph{Compatibility of {Kisin} modules for different uniformizers},
  Journal f\"ur die reine und angewandte Mathematik (2013).

\bibitem{Morrow-TsujiGeneralizedrep}
Matthew Morrow and Takeshi Tsuji, \emph{Generalised representations as
  {$q$}-connections in integral {$p$}-adic hodge theory}, 2020.

\bibitem{Ozeki}
Yoshiyasu Ozeki, \emph{Lattices in crystalline representations and {K}isin
  modules associated with iterate extensions}, Doc. Math. \textbf{23} (2018),
  497--541. \MR{3846051}

\bibitem{ScholzeWeinstein}
Peter Scholze and Jared Weinstein, \emph{Berkeley lectures on {$p$}-adic
  geometry}, Annals of Math Studies 207.

\bibitem{stacks-project}
The {Stacks Project Authors}, \emph{{Stacks Project}},
  \url{http://stacks.math.columbia.edu}, 2021.

\end{thebibliography}

\end{document}